\documentclass[10pt]{article}
\usepackage{amsmath, amssymb}
\usepackage{amsthm} 
\usepackage[all]{xy}

\def\Sch{\text{\rm Sch}}
\def\F{\cal F}
\def\L{{\cal L}}

\def\bgn{\begin}

\def\CL{\text{\rm CL}}

\def\J{{\cal J}}
\def\L{{\cal L}}

\def\1{{[1]}}
\def\2{{[2]}}
\def\3{{[3]}}
\def\({\left(}
\def\){\right)}
\def\s-circ{\,{\scriptstyle{\circ}}\,}
\def\<<{<\negthinspace \negthinspace<}
\def\Ad{\text{\rm Ad}}
\def\ad{\text{\rm ad}}

\def\even{\text{\rm even}}
\def\bgn{\begin}

\def\endaln{\end{align}}
\def\cou{\ss\text{\rm co}}

\def\<{<\negthinspace \negthinspace <}

\def\t{\theta}

\def\({\left(}
\def\){\right)}
\def\Im{\text{\rm Im}}
\def\Re{\text{\rm Re}}
\def\[{\big[\neg\big[}
\def\]{\big]\neg\big]}
\def\al{\al}

\def\M{{\cal M}}

\def\tr{\text{\rm tr}}
\def\a{\alpha}
\def\b{\beta}

\def\e{\varepsilon}

\def\Gam{\Gamma}
\def\k{\kappa}
\def\del{\delta}
\def\lam{\lambda}
\def\ome{\omega}

\def\sig{\sigma}

\def\A{{\cal A} }

\def\G{G_2}
\def\D{\Bbb D}

\def\R{\Bbb R}
\def\C{\Bbb C}

\def\M{\frak M}

\def\w{\wedge}
\def\({\left(}
\def\){\right)}
\def\G{G_2}

\def\neg{\negthinspace}
\def\h{\hat}

\def\wtil{\widetilde}

\def\ol{\overline}
\def\pa{\partial}

\def\ran{\rangle} 
\def\lan{\langle}

\def\ss{\scriptscriptstyle}

\def\arrow{\longrightarrow}

\def\:{\, :\,}
\def\CL{\text{\rm CL}}
\def\TT{T\oplus T^*}
\def\complex{generalized complex }
\def\K\"ahler{generalized K\"ahler}
\def\vol{\text{\rm vol}}

\def\10{\displaystyle L^{10}}
\def\2{\displaystyle L^2}
\def\c0{\displaystyle C^0}

\def\10{\displaystyle L^{10}}
\def\2{\displaystyle L^2}
\def\del{\delta}
\def\del2{\displaystyle L^2_{0,\delta}}
\def\c0{\displaystyle C^0}

\def\del{\delta}

\def\K{{\cal K}}
\def\M-A{\text{\rm Monge-Amp\`ere}}

\def\M-A{\text{\rm Monge-Amp\`ere}}
\def\[{\big[\,}
\def\]{\,\big]}
\def\End{\text{\rm End\,}}

\def\id{\text{\rm id}}

\theoremstyle{plain} 
\newtheorem{theorem}{\indent\sc Theorem}[section] %

\newtheorem{proposition}[theorem]{\indent\sc Proposition}

\theoremstyle{definition}

\makeatletter
\def\even{\text{\rm even}}
\def\odd{\text{\rm odd}}
\def\GK{\text{\rm generalized K\"ahler }}
\def\L{\text{{\cal L}}}
\def\TT{TM\oplus T^*M}

\def\M{{\cal M}}

\def\ol{\overline}
\def\part{\partial}
\def\Con{\wtil{\cal{M}}}
\def\G{{\cal G}}
\def\GC{{\cal{GM}}}
\def\GCon{\wtil{\cal GM}}
\def\tt{\ss T\oplus T^*}
\def\D{{\cal D}}
\def\F{{\cal F}}
\def\L{{\cal L}}
\def\Lam{\Lambda}
\def\Herm{\text{\rm Herm}}
\def\prim{\text{\rm prim}}
\def\ad{\text{\rm ad}}
\title{Moduli spaces of Einstein-Hermitian generalized connections over generalized K\"ahler manifolds of symplectic type
}%
\author{Ryushi Goto}
\date{} 
\begin{document}
\maketitle

\begin{abstract}
From a view point of the moment map, we shall introduce the notion of Einstein-Hermitian generalized connections over a generalized K\"ahler manifold of symplectic type.
 We show that
 moduli spaces of Einstein-Hermitian generalized connections arise as the K\"ahler quotients. 
The deformation complex of Einstein-Hermitian generalized connections is
an elliptic complex and it turns out that the smooth part of the moduli space is a finite dimensional K\"ahler manifold. 
The canonical line bundle over a generalized K\"ahler manifold of symplectic type has the canonical generalized connection and its curvature coincides with "the scalar curvature as the moment map" which is defined in the previous paper \cite{Goto_2016}. 
K\"ahler-Ricci solitons provide examples of Einstein-Hermitian generalized connections and Einstein Hermitian co-Higgs bundles are also discussed. 
\end{abstract}
\tableofcontents

\numberwithin{equation}{section}
\section{Introduction}
Behind the K\"ahler geometry, a framework of symplectic geometry appears and plays an important role.
Fujiki and Donaldson showed that the scalar curvature of K\"ahler manifolds arises as the moment map \cite{Fu_1990}
, \cite{Do_1997}. 
 The mean curvature of Hermitian connections over a compact K\"ahler manifold is also regarded as the moment map. Atiyah-Bott constructed a K\"ahler structure on the moduli spaces of stable vector bundles over compact Riemannian surfaces \cite{AB_1983}.
S. Kobayashi obtained moduli spaces of Einstein-Hermitian vector bundles over compact K\"ahler manifolds which arise as K\"ahler quotients \cite{Ko_1987}. 
In the previous paper \cite{Goto_2016}, the author obtained the scalar curvature as moment map over generalized K\"ahler manifolds.
In this paper we will pursue an analogue of the moment map framework for generalized Hermitian connections over a generalized K\"ahler manifold and construct moduli spaces of Einstein-Hermitian generalized connections as the K\"ahler quotients. 
 A generalized connection $\D^A$ over a complex vector bundle $E$ is a differential operator from $\Gam(E)$ to $\Gam(E\otimes(\TT))$ which satisfies  
 $$
 \D^\A(f s) = s\otimes df  + f\D^\A s, \qquad  s\in \Gam(E),\,\,\,f\in C^\infty(M),
 $$
A generalized connection consists of an ordinary
  connection $D^A$ and $\End(E)$-valued vector field $V$, i.e., 
  $\D^\A=D^A+V.$
While the notion of the curvature of generalized connections has been expected, a precise
formulation of the curvature satisfying the moment map framework
has not yet appeared in the literature.
In this paper, a $d$-closed non-degenerate, pure spinor $\psi$ on $M$ plays a central role to obtain the curvature of a generalized connection.
For simplicity we mainly consider a $d$-closed non-degenerate, pure spinor $\psi$ on $M$ ${}^{\dagger}$ which is given by 
the exponential $e^{b+\sqrt{-1}\ome}$, where 
$\ome$ is a real symplectic form and 
$b$ is a $d$-closed real $2$-form, that is, 
\footnote {${}^{\dagger}$We can not give a definition of the curvature of generalized connections in full generality. However, our construction of the curvature and main theorems apply to the class of generalized K\"ahler manifolds of symplectic type which constructed from Poisson deformations \cite{Goto_2010}.}
$$
\psi=1+ (b+\sqrt{-1}\ome) +\frac1{2!}(b+\sqrt{-1}\ome)^2+\cdots+\frac1{n!}(b+\sqrt{-1}\ome)^n.
$$
A $d$-closed, nondegenerate, pure spinor $\psi$ induces the \complex structure $\J_\psi$.
{\it A generalized K\"ahler structure of symplectic type} is a generalized K\"ahler structure which is a pair of generalized complex  structures 
$(\J, \J_\psi)$, where $\J$ is an arbitrary \complex structure (see Section 2 for details).
A vector field acts on the differential form $\psi$ by the interior product and 
a  differential form also acts on $\psi$ by exterior product. 
Both actions of multi-vector fields and differential forms give the spin representation of the Clifford algebra.
Then the spin representation together with the trivial action of $\End(E)$ 
yields the action of $\End(E)\otimes \w^\bullet(\TT)$ on $\psi$.
Then we define a curvature $\F_\A(\psi)$ of a generalized connection $\D^\A=D^A+V$ to be an $\End(E)$-valued differential form which is given by
$$
\F_\A(\psi) =F_A\cdot\psi+ d^{D_A}(V\cdot\psi)+\frac12 [V, V]\cdot \psi,
$$
where the curvature $F_A$ of the ordinary connection $D^A$ is an $\End(E)$-valued $2$-form which acts on $\psi$ and 
the action of $V$ on $\psi$  gives an $\End(E)$-valued form $V\cdot\psi$ and 
the covariant derivative $d^{D_A}$ of the ordinary connection  $D^A$ also gives an $\End(E)$-valued form  $d^{D_A}(V\cdot\psi)$. 
The commutator $[V,V]$ is an $\End(E)$-valued $2$-vector acting on $\psi.$
We also have the notion of Hermitian generalized connection of a Hermitian vector bundle $(E, h)$ over a generalized K\"ahler manifold $(M,\J, \J_\psi)$. We denote by $\GCon(E,h)$ the affine space of generalized Hermitian connections of $(E,h)$. Then it turns out that $\GCon(E,h)$ has a natural K\"ahler structure  whose symplectic form $\ome_{\GCon}$ is given by 
\bgn{equation}\label{sym}
\ome_{\GCon}(\dot{\A}_1, \dot\A_2)=-\int_M \tr \lan \J_\psi\dot\A_1, \dot\A_2\ran_{\tt} \vol_M,
\end{equation}
where $\dot{\A}_1, \dot{\A}_2\in \Gam(u(E)\otimes (\TT))$ are tangents of the space $\GCon(E,h)$ (see Section 6 for details).
Then the unitary gauge group $U(E,h)$ acts on $\GCon(E,h)$ preserving the symplectic structure $\ome_{\GCon}.$
Then we have \par\smallskip
{\indent\sc Proposition 6.1. }{\it
Let  $\mu:\GCon(E, h)\to u(E)^*$ be the moment map on $\GCon(E, h)$ for the action of 
$U(E, h)$. Then the moment map $\mu$ is given by
$$
\lan\mu(\A),\xi\ran=\int_M \Im \,\,i^{-n}\tr \lan\xi\psi,\,\, \F_\A(\ol\psi)\ran_s,
$$
where $\mu(\A)\in u(E)^*$ denotes the value of $\mu$ at a generalized connection $\D^\A\in \GCon(E,h)$ and the coupling $\mu(\A)$ and
$\xi \in u(E)$ is denoted by $\lan\mu(\A),\xi\ran$ and $\Im$ is the imaginary part of $i^{-n}\tr \lan\xi\psi,\,\, \F_\A(\ol\psi)\ran_s.$}\\
\smallskip
\par
We define a mean curvature $\K_\A(\psi)$ as 
$$
\K_\A(\psi)=\pi_{U^{-n}}^{\Herm}\F_\A(\psi),
$$
where $\pi_{U^{-n}}^{\Herm}$ denotes the projection to the component $\Herm(E, h)\otimes\psi$, where 
$\Herm(E,h)$ denotes Hermitian endmorphisms of $E$.
Then the Einstein-Hermitian condition is given by 
\par\smallskip
{\noindent\sc Definition 3.3. }A generalized Hermitian connection $\D^\A$ is {\it Einstein-Hermitian} if $\D^\A$
satisfies the following: 
$$
\K_\A(\psi)=\lam id_E, \quad \text{\rm for a real constant }\lam
$$
\par\smallskip
Under the identification $u(E)^*$ with $\Herm(E,h)$, the moment map $\mu$ is given by 
$$
\mu(\A) =\K_\A(\psi):=\pi^{\Herm}_{U^{-n}}\F_\A(\psi)
$$
A generalized Hermitian connection $\D^\A$ is decomposed as $\D^\A=\D^{1,0}+\ol\pa^\A_\J,$
where $\ol\pa^\A_\J : E\to E\otimes \ol{\L_\J}$. 
If $\ol\pa^\A_J$ satisfies $\ol\pa^\A_\J\circ \ol\pa^\A_\J=0$, then $\ol\pa^\A_\J$ is called {\it a generalized holomorphic structure} of $E$.
This is an analogue of the notion of Hermitian connections associated with holomorphic vector bundles which is introduced and developed by \cite{Gua_2010}, \cite{Hi_2011-1}.
We define $\GCon(E, h)^{hol}$ to be a subset of $\GCon(E, h)$ consisting of $\D^\A$ such that 
$\ol\pa^\A_\J\circ \ol\pa^\A_\J=0$, that is, 
$$
\GCon(E, h)^{hol}:=\{\, \D_\A\in \GCon(E,h)\, |\, \ol\pa^\A_\J\circ \ol\pa^\A_\J=0\, \}
$$
Then it turns out that 
$\GCon(E, h)^{hol}$ is a complex submanifold of $\GCon(E, h)$ which inherits the K\"ahler structure.
Then 
we define the set of Einstein-Hermitian generalized connections $ {\GCon(E, h)}^{hol}_{EH}$  by 
$$
 {\GCon(E, h)}^{hol}_{EH}=
\{\, \D_\A\in \GCon(E, h)^{hol}\, |\, \K_\A(\psi)=\lam \id_E\, \}
$$
The moduli space of Einstein-Hermitian generalized connections ${{\cal GM}(E, h)}^{hol}_{EH}$ is defined to be the quotient 
of ${{\cal GM}(E, h)}^{hol}_{EH}$ by the action of the gauge group $U(E, h)$, that is, 
$$
{{\cal GM}(E, h)}^{hol}_{EH}= {\GCon(E, h)}^{hol}_{EH}/\G(E, h)
$$
\medskip
Then our main theorem is the following.
\par\smallskip{\noindent\sc Theorem 6.8. }{\it 
The moduli space ${{\cal GM}(E, h)}^{hol}_{EH}$ arises as 
a K\"ahler quotient and the smooth part of ${{\cal GM}(E, h)}^{hol}_{EH}$ inherits a K\"ahler structure.}
\par\medskip
In Section 2, we will give a brief review of generalized complex structures and generalized K\"ahler structures.
In Section 3, the curvature $\F_\A(\psi)$ of a generalized connection $\D^\A$  of a vector bundle $E$ on a symplectic manifold $(M,\ome)$ is introduced${}^{\ddagger}$.\footnote{${}^{\ddagger}$Our definition of the curvature of generalized connection can be applied to the cases of 
arbitrary $d$-closed nondegenerate, pure, spinors.} The differential form $\tr\F_\A(\psi)$ is a representative of 
the first Chern class of 
the vector bundle $E$ coupled with the class $[\psi]$.
In Section 4, we explain that there is  a unique generalized Hermitian connection associated with a generalized holomorphic structure 
which is called the canonical generalized connection.
In Section 5, we discuss the canonical generalized connection of the canonical line bundle $K_\J$ which is given by a generalized complex structure $\J$ over a generalized K\"ahler manifold of symplectic type $(\J, \J_\psi). $ Then $K_\J$ is a Hermitian line bundle which admits a generalized holomorphic structure. We show that that the curvature of the canonical generalized connection of $K_\J$ coincides with "the scalar curvature as the moment map" 
which is shown in the previous paper \cite{Goto_2016}. We also show that a K\"ahler-Ricci soliton provides an Einstein-Hermitian generalized connection over a generalized K\"ahler manifold. 
In Section 6, we  show that $\K_\A(\psi)$ constructed from the curvature $\F_\A(\psi)$ is the moment map $\mu$ and 
the moduli space ${{\cal{GM}}(E, h)}^{hol}_{EH}$ is constructed. 
In Section 7, we show that the deformation complex of  Einstein-Hermitian generalized connections is an elliptic complex. 
In Section 8, we will discuss the Einstein-Hermitian condition of co-Higgs bundles.
\section{Generalized complex structures and generalized K\"ahler structures}
\subsection{Generalized complex structures and nondegenerate, pure spinors}
Let $M$ be a differentiable manifold of real dimension $2n$.
The bilinear form $\lan\,\,,\,\,\ran_{\tt}$ on 
the direct sum $T_M \oplus T^*_M$ over a differentiable manifold $M$ of dim$=2n$ is defined by 
$$\lan v+\xi, u+\eta \ran_{\tt}=\frac12\(\xi(u)+\eta(v)\),\quad  u, v\in T_M, \xi, \eta\in T^*_M .$$
Let SO$(\TT)$ be the fibre bundle over $M$ with fibre SO$(2n, 2n)$ which is 
a subbundle of End$(\TT)$  preserving the bilinear form $\lan\,,\,\ran_s$ 
 An almost \complex structure $\J$ is a section of SO$(\TT)$ satisfying $\J^2=-\id.$ Then as in the case of almost complex structures, an almost \complex structure $\J$ yields the eigenspace decomposition :
$(T_M \oplus T^*_M)^\C =\L_\J \oplus \ol \L_\J$, where 
$\L_\J$ is $-\sqrt{-1}$-eigenspaces and  $\ol{\L}_\J$ is the complex conjugate of $\L_\J$. 
The Courant bracket of $\TT$ is defined by 
$$
 [u+\xi, v+\eta]_{\cou}=[u,v]+{\cal L}_u\eta-{\cal L}_v\xi-\frac12(di_u\eta-di_v\xi),
 $$
 where $u, v\in TM$ and $\xi, \eta$ is $T^*M$.
If $\L_\J$ is involutive with respect to the Courant bracket, then $\J$ is a generalized complex structure, that is, $[e_1, e_2]_{\cou}\in \Gam(\L_\J)$  for any two elements 
 $e_1=u+\xi,\,\, e_2=v+\eta\in \Gam(\L_\J)$.
Let $\CL(T_M \oplus T^*_M)$ be the Clifford algebra bundle which is 
a fibre bundle with fibre the Clifford algebra $\CL(2n, 2n)$ with respect to $\lan\,,\,\ran_{\tt}$ on $M$.
Then a vector $v$ acts on the space of differential forms $\oplus_{p=0}^{2n}\w^pT^*M$ by 
the interior product $i_v$ and a $1$-form acts on $\oplus_{p=0}^{2n}\w^pT^*M$ by the exterior product $\t\w$, respectively.
Then the space of differential forms gives a representation of the Clifford algebra $\CL(\TT)$ which is 
the spin representation of $\CL(\TT)$. 
Thus
the spin representation of the Clifford algebra arises as the space of differential forms $$\w^\bullet T^*_M=\oplus_p\w^pT^*_M=\w^{\even}T^*_M\oplus\w^{\odd}T^*_M.$$ 
The inner product $\lan\,,\,\ran_s$ of the spin representation is given by 
$$
\lan \a, \,\,\,\b\ran_s:=(\a\w\sig\b)_{[2n]},
$$
where $(\a\w\sig\b)_{[2n]}$ is the component of degree $2n$ of $\a\w\sig\b\in\oplus_p \w^pT^*M$ and 
$\sig$ denotes the Clifford involution which is given by 
$$
\sig\b =\bgn{cases}&+\b\qquad \deg\b \equiv 0, 1\,\,\mod 4 \\ 
&-\b\qquad \deg\b\equiv 2,3\,\,\mod 4\end{cases}
$$
We define $\ker\Phi:=\{ e\in (T_M\oplus T^*_M)^\C\, |\, e\cdot\Phi=0\, \}$ for a differential form $\Phi
\in \w^{\even/\odd}T^*_M.$
If $\ker\Phi$ is maximal isotropic, i.e., $\dim_\C\ker\Phi=2n$, then $\Phi$ is called {\it a pure spinor} of even/odd type.
A pure spinor $\Phi$ is {\it nondegenerate} if $\ker\Phi\cap\ol{\ker\Phi}=\{0\}$, i.e., 
$(T_M\oplus T^*_M)^\C=\ker\Phi\oplus\ol{\ker\Phi}$.
Then a nondegenerate, pure spinor $\Phi\in \w^\bullet T^*_M$ gives an almost generalized complex structure $\J_{\Phi}$ which satisfies 
$$
\J_\Phi e =
\bgn{cases}
&-\sqrt{-1}e, \quad e\in \ker\Phi\\
&+\sqrt{-1}e, \quad e\in \ol{\ker\Phi}
\end{cases}
$$
Conversely, an almost \complex structure $\J$ locally arises as $\J_\Phi$ for a nondegenerate, pure spinor $\Phi$ which is unique up to multiplication by
non-zero functions.  Thus an almost \complex structure yields the canonical line bundle $K_{\J}:=\C\lan \Phi\ran$ which is a complex line bundle locally generated by a nondegenerate, pure spinor $\Phi$ satisfying 
$\J=\J_\Phi$.
An \complex structure 
$\J_\Phi$ is integrable if and only if $d\Phi=\eta\cdot\Phi$ for a section $\eta\in T_M\oplus T^*_M$. 
The {\it type number} of $\J=\J_\Phi$ is defined as the minimal degree of the differential form $\Phi$. Note that type number Type $\J$ is a function on a manifold which is not a constant in general.
\bgn{example}
Let $J$ be a complex structure on a manifold $M$ and $\J^*$ the complex structure on the dual  bundle $T^*M$ which is given by $J^*\xi(v)=\xi (Jv)$ for $v\in TM$ and $\xi\in T^*M$.
Then a \complex structure $\J_J$ is given by the following matrix
$$\J_J=\bgn{pmatrix}J&0\\0&-J^*
\end{pmatrix},$$
Then the canonical line bundle is the ordinary one which is generated by complex forms of type $(n,0)$.
Thus we have  Type $\J_J =n.$
\end{example}
\bgn{example}
Let $\ome$ be a symplectic structure on $M$ and $\h\ome$ the isomorphism from $TM$ to $T^*M$ given by $\h\ome(v):=i_v\ome$. We denote by $\h\ome^{-1}$ the inverse map from $T^*M$ to $TM$.
Then a \complex structure $\J_\psi$ is given by the following
$$\J_\psi=\bgn{pmatrix}0&-\h\ome^{-1}\\
\h\ome&0
\end{pmatrix},\quad\text{\rm Type $\J_\psi =0$}$$
Then the canonical line bundle is given by the differential form $\psi=e^{\sqrt{-1}\ome}$. 
Thus Type $\J_\psi=0.$
\end{example}
\bgn{example}[$b$-field action]
A $d$-closed $2$-form $b$ acts on a \complex structure by the adjoint action of Spin group $e^b$ which provides
a \complex structure $\Ad_{e^b}\J=e^b\circ \J\circ e^{-b}$. 
\end{example}
\bgn{example}[Poisson deformations]
Let $\b$ be a holomorphic Poisson structure on a complex manifold. Then the adjoint action of Spin group $e^\b$ gives deformations of new \complex structures by 
$\J_{\b t}:=\Ad_{\b^{Re} t}\J_J$.  Then Type ${\J_{\b t}}_x=n-2$ rank of $\b_x$ at $x\in M$,
which is called the Jumping phenomena of type number.
\end{example}
Let $(M, \J)$ be a generalized complex manifold and $\ol \L_\J$ the eigenspace of eigenvalue $\sqrt{-1}$.
Then we have the Lie algebroid complex $\w^\bullet\ol{\L}_\J$:
$$
0\arrow\w^0\ol \L_\J\overset{\ol\pa_\J}\arrow\w^1\ol \L_\J\overset{\ol\pa_\J}\arrow\w^2\ol \L_\J\overset{\ol\pa_\J}\arrow\w^3\ol \L_\J\arrow\cdots
$$
The Lie algebrid complex is the deformation complex of \complex structures. 
In fact, $\e\in \w^2\ol \L_\J$ gives deformed isotropic subbundle 
$E_\e:=\{ e+[\e, e]\, |\, e\in \L_\J\}$. 
Then $E_\e$ yields deformations of \complex structures if and only if $\e$ satisfies Generalized Mauer-Cartan equation
$$
\ol{\pa}_\J\e+\frac12[\e, \e]_{\Sch}=0,
$$
where $[\e, \e]_{\Sch}$ denotes the Schouten bracket. 
The Kuranishi space of generalized complex structures is constructed.
Then the second cohomology group $H^2(\w^\bullet\ol \L_\J)$ of the Lie algebraic complex gives the infinitesimal deformations of \complex structures and the third one 
$H^3(\w^\bullet\ol \L_\J)$ is the obstruction space to deformations of \complex structures.
Let $\{e_i\}_{i=1}^n$ be a local basis of $\L_\J$ for an almost \complex structure $\J$, 
where $\lan e_i, \ol e_j\ran_{\tt}=\del_{i,j}$.
The the almost \complex structure $\J$ is written as an element of Clifford algebra,
$$
\J=\frac{\sqrt{-1}}2\sum_i e_i\cdot\ol {e}_i,
$$
where $\J$ acts on $\TT$ by the adjoint action $[\J, \,]$. 
Thus we have $[\J, e_i]=-\sqrt{-1}e_i$ and $[\J, \ol e_i]=\sqrt{-1}e_i$.
An almost \complex structure $\J$ acts on differential forms by the Spin representation which gives the decomposition into eigenspaces:
\bgn{equation}
\w^\bullet T^*_M=U^{-n}\oplus U^{-n+1}\oplus\cdots U^{n},
\end{equation}
where $U^{i}(=U^i_\J)$ denotes the $i$-eigenspace. 
\subsection{Generalized K\"ahler structures}
\bgn{definition}
{\it A generalized K\"ahler structure} is a pair $(\J_1, \J_2)$ consisting of two commuting \complex structures 
$\J_1, \J_2$ such that $\h G:=-\J_1\circ\J_2=-\J_2\circ \J_1$ gives a positive definite symmetric form 
$G:=\lan \h G\,\,,  \,\,\ran$ on $T_M\oplus T_M^*$, 
We call $G$ {\it a generalized metric}.
{\it A \GK structure of symplectic type} is a \GK structure $(\J_1, \J_2)$ such that  $\J_2$ is the \complex structure $\J_\psi$ which is  induced from  a $d$-closed, nondegenerate, pure spinor $\psi:=e^{b+\sqrt{-1}\ome}.$
\end{definition}
Each $\J_i$ gives the decomposition $(\TT)^\C=\L_{\J_i}\oplus\ol \L_{\J_i}$ for $i=1,2$.
Since $\J_1$ and $\J_2$ are commutative, we have the simultaneous eigenspace decomposition 
$$
(\TT)^\C=(\L_{\J_1}\cap \L_{\J_2})\oplus (\ol \L_{\J_1}\cap \ol \L_{\J_2})\oplus (\L_{\J_1}\cap \ol \L_{\J_2})\oplus
(\ol \L_{\J_1}\cap \L_{\J_2}).
$$
Since $\h G^2=+\id$,
The generalized metric $\h G$ also gives the eigenspace decomposition: $\TT=C_+\oplus C_-$, 
where $C_\pm$ denote the eigenspaces of $\h G$ of eigenvalues $\pm1$. 
We denote by $\L_{\J_1}^\pm$ the intersection $\L_{\J_1}\cap C^\C_\pm$. 
Then it follows 
\bgn{align*}
&\L_{\J_1}\cap \L_{\J_2}=\L_{\J_1}^+,  \quad \ol \L_{\J_1}\cap \ol \L_{\J_2}=\ol \L_{\J_1}^+\\
&\L_{\J_1}\cap \ol \L_{\J_2}=\L_{\J_1}^-,\quad \ol \L_{\J_1}\cap \L_{\J_2}=\ol \L_{\J_1}^-
\end{align*}
\bgn{example}
Let $X=(M, J,\ome)$ be a K\"ahler manifold. Then the pair $(\J_J, \J_\psi)$ is a generalized K\"ahler where 
$\psi=\exp(\sqrt{-1}\ome)$. 
\end{example}
\bgn{example}
Let $(\J_1, \J_2)$ be a generalized K\"ahler structure. 
Then the action of $b$-fields gives a generalized K\"ahler structure 
$(\Ad_{e^b}\J_1, \Ad_{e^b}\J_2)$ for a $d$-closed $2$-form $b.$
\end{example}
\subsection{The stability theorem of generalized K\"ahler manifolds}
It is known that the stability theorem of ordinary K\"ahler manifolds holds
\bgn{theorem}[Kodaira-Spencer]
Let $X=(M,J)$ be a compact K\"ahler manifold and $X_t$ small deformations of $X=X_0$ as complex manifolds.
Then $X_t$ inherits a K\"ahler structure. 
\end{theorem}
The following stability theorem of generalized K\"ahler structures provides many interesting examples of generalized K\"ahler manifolds.
\bgn{theorem}{\rm \cite{Goto_2010}}
Let $X=(M,J,\ome)$ be a compact K\"ahler manifold and $(\J_J, \J_\psi)$ the induced \K\"ahler structure, 
where $\psi=e^{\sqrt{-1}\ome}$. 
If there are analytic deformations $\{\J_t\}$ of $\J_0=\J_J$ as \complex structures, then there are deformations of $d$-closed nondegenerate, pure spinors $\{\psi_t\}$ such that 
pairs $(\J_t, \J_{\psi_t})$ are \K\"ahler structures, where $\psi_0=\psi$
\end{theorem}
\section{Einstein-Hermitian generalized connections}
\subsection{Generalized connections over symplectic manifolds}
Let $M$ be a compact real manifold of dimension $2n$ and $\ome$ a real symplectic structure on $M$. 
We denote by $\psi$ the exponential of $b+\sqrt{-1}\ome$, that is, 
$$
\psi:=e^{b+\sqrt{-1}}=1+(b+\sqrt{-1}\ome)+\frac1{2!}(b+\sqrt{-1}\ome)^2+\cdot+\frac1{n!}(b+\sqrt{-1}\ome)^n,
$$
where $b$ denotes a $d$-closed $2$-form on $M.$
Then $\psi$ induces a generalized complex structure $\J_\psi$ which gives a decomposition 
$$
(\TT)^{\Bbb C}=\L_{\J_\psi}\oplus \ol {\L_{\J_\psi}}
$$
$\TT$ acts on differential forms by the spin representation which is given by the interior product and the exterior product of $\TT$ on differential forms. 
Then we have a decomposition of differential forms on $M,$
$$
\oplus_{i=0}^{2n}\w^iT^*_M=\oplus_{j=0}^{2n} U_{\J_\psi}^{-n+j}, 
$$
where $U_{\J_\psi}^{-n}$ is a complex line bundle generated by $\psi$  and 
$U_{\J_\psi}^{-n+i}$ is constructed by the spin action of $\w^i\ol{\L_{\J_\psi}}$ on $U^{-n}_{\J_\psi}.$
Let $E\to M $ be a complex vector bundle of rank $r$ over $M$ and $\Gam(E)$ a set of smooth sections of $E$.
We denote by $\Gam(E\otimes(\TT)^\C)$ the set of smooth sections of $E\otimes(\TT)^\C$. 
{\it A generalized connection}
$\D^\A$
is a map from $\Gam(E)$ to $\Gam(E\otimes(\TT)^\C)$ such that 
$$
\D^\A(fs)=s\otimes df+f\D^\A(s),\qquad \text{\rm for } s\in\Gam(E),\,\, f\in C^\infty(M).
$$
Let $h$ be a Hermitian metric of $E$. {\it A generalized Hermitian connection } is a generalized connection $\D^\A$ satisfying 
$$
dh(s, s')=h(\D^\A s, s')+h(s, \D^\A s'), \qquad \text{\rm for } s, s'\in \Gam(E).
$$
We denote by $u(E,h)(:=u(E))$ the set of skew-symmetric endomorphisms of $E$ with respect to $h.$ Then $\End(E)$ is decomposed as
\bgn{equation}\label{eq:End(E)}
\End(E)=u(E,h)\oplus \Herm(E,h),
\end{equation}
where $\Herm(E,h)$ denotes the set of Hermitian endmorphisms of $E.$
Let $\{U_\a\}$ be an open covering of $M$ which gives local trivializations of $E$. 
We take  $s_\a:=(s_{\a,1}, \cdots, s_{\a,r})$ as  a local unitary frame of $E$ over  $U_\a$. 
The set of transition functions is denoted by $\{g_{\a\b}\}$.
Then given a connection $\D^\A$ is written as 
$$\D^\A(s_{\a,p})=\sum_{q=1}^r s_{\a,q}\A_{p,\a}^q,
$$
where $\A_\a:=\{\A_{p, \a}^q\}$ is a section of $u(E)\otimes_\R(\TT)|_{U_\a}$ which is denoted as 
$$
\A_\a=\sum_i \A_\a^i e_i,
$$
where $e_i\in (\TT)$ and $\A_\a^i\in u(E).$
Note that each $e_i$ is a real element of $\TT$.
Then $\A_\a$ is also decomposed into  
$$\A_{p,\a}^q=A_{p,\a}^q+V_{p,\a}^q.$$  
Then it turns out that$A_\a:=\{A_{p,\a}^q\}$ is a connection form and 
$V_\a:=\{V_{p,\a}^q\}$ gives a section of $u(E)\otimes T_M$
In fact, 
by using local trivializations $s_\a$, given a generalized connection $\D^\A$ is written as 
$$
\D_\A= d_\a+  \A_\a,
$$
Since 
$\D_\A :\Gam (E)\to \Gam (E\otimes (\TT))$ is globally defined as a differential operator, we have 
\bgn{align}
A_\a=&-(dg_{\a\b})g^{-1}_{\a\b}+g_{\a\b}A_\b g^{-1}_{\a\b}\\
V_\a=&g_{\a\b}V_\b g^{-1}_{\a\b}
\end{align}
Thus it follows that $A_\a$ is a connection form of a connection $D^A$ and  $V_\a$ is a $u(E)$-valued vector field. 
As shown in the ordinary connections, a generalized connection $\D_\A $ is extended to an operator　$\Gam(\End(E))\to \Gam(\End(E)\otimes(\TT))$
by the following: 
$$
(\D_\A \xi) s:=\D_\A(\xi s)-\xi(\D_\A s), \qquad \text{\rm for }\xi\in \Gam(\End(E)), \,\, s\in\Gam(E).
$$
Then the extended operator $\D^\A$  is also written as follows in terms of  $\A_\a=\sum_i\A_{\a,i}e_i$　
\bgn{align}
\D_\A\xi=d\xi+\sum_i[\A_{\a,i}, \xi]e_i.
\end{align}
 Each $e_i\in \TT$ acts on $\psi$ which is denoted by $e_i\cdot\psi\in U^{-n+1}$ and then each element $\xi\otimes e_i$ of 
$\End(E)\otimes(\TT)$ acts on $\psi$ denoted by $\xi\otimes (e_i\cdot\psi).$ 
Then the action on $\psi$ gives the map 
$\End(E)\otimes(\TT)\to \End(E)\otimes U^{-n+1}$. 
Thus we obtain $\D^\A\xi \cdot\psi \in \End(E)\otimes U^{-n+1}$ for $\xi\in \End(E)$ and a generalized connection $\D^\A.$
Thus 
 an operator $d^{\D_\A}$ is defined by  
 \bgn{align}
 d^{\D_\A}: \Gam(\End&(E))\to\Gam(\End(E)\otimes U^{-n+1})\\
& \xi \mapsto \D^\A\xi\cdot\psi 
 \end{align}
We also extend $d^{\D_\A}$ to be an operator 
$$d^{\D_\A}:\Gam(\End(E)\otimes U^{-n+1})\to \Gam(\End(E)\otimes( U^{-n}\oplus U^{-n+2}))$$
by setting :
$$
d^{\D_\A}(\xi\otimes e_i\cdot\psi)=(\D^\A\xi)\cdot e_i\cdot\psi+\xi\otimes d(e_i\cdot\psi)
$$
Let $a=\sum_i a_i e_i$ be a section of $\End(E)\otimes(\TT)$, where 
$a_i\in \End(E)$ and $e_i\in \TT$.
Then the extended operator $d^{\D_\A}$ is written in the following form:
\bgn{align}\label{dDA(acdotpsi)}
d^{\D_\A}(a\cdot\psi)=&\sum_i d^{\D_\A}(a_ie_i\cdot\psi)=\sum_i (d^{\D_\A}a_i)e_i\cdot\psi+
a_id(e_i\cdot\psi)\notag\\
=&\sum_{i,j}da_i\cdot e_i\cdot \psi+a_id(e_i\cdot\psi)+[\A_{\a, j}, \, a_i]e_j\cdot e_i\cdot\psi\notag\\
=&d(a \cdot\psi)+[\A\, \cdot a]\cdot\psi
\end{align}
where  
$[\A\, \cdot a]:=\sum_{i,j}[A_{\a, j}, \, a_i]e_i\cdot e_j.$
By using the local trivialization and the decomposition $\A_\a=A_\a+V_\a$, the operator $d^{\D_\A}$ is described as the following:
$$
d^{\D_\A} (a\cdot\psi) =(d_\a+A_\a+ V_\a) a\cdot\psi =d(a\cdot\psi)+ [V_\a\cdot a]\cdot\psi+ [A_\a\cdot a]\cdot\psi,
$$
where $[\,\cdot\,]$ is the product of 
$\End(E)\otimes \CL$ which is defined as 
$$[V_\a\cdot a]=\sum_{i,j}[V_{\a, i}, a_j]e_i\cdot e_j$$
$$[A_\a\cdot a]=\sum_{i,j}[A_{\a, i}, a_j]e_i\cdot e_j$$

Note that  our new product $[\, \cdot\, ]$ is the combinations of 
 the bracket of 
Lie algebra $\End(E)$ and the Clifford multiplications of Clifford algebra
$\CL$ gives , so that is, 
$$
[(A\otimes s)\cdot (A'\otimes s')] =[A, A']\otimes s\cdot s'\in \End(E)\otimes\CL
$$
for $ A \otimes e, A'\otimes e'$, where $A, A'\in \End(E), 
s, s'\in \CL$.
\subsection{Curvature of generalized connections}
Let $\D^\A$ be a generalized Hermitian connection of a Hermitian vector bundle $(E,h)$ over a manifold $M$ which consists an ordinary Hermitian connection $D^A$ and 
a section $V\in \Gam(u(E)\otimes T_M)$. 
Then the ordinary curvature $F_A$ is a section $\End(E)$ valued $2$-from which acts on $\psi$ by the spin
representation to obtain $F_A\cdot\psi \in \End(E)\otimes_\R (U^{-n}\oplus U^{-n+2}).$
The ordinary connection is extended to be an operator as before
$$
d^A: \End(E)\otimes_\R U^{-n+1}\to \End(E)\otimes_\R (U^{-n}\oplus U^{-n+2})
$$
By applying $d^A$ to $V\cdot\psi \in \End(E)\otimes_\R U^{-n+1}$, we have $d^A(V\cdot\psi)\in \End(E)\otimes_\R (U^{-n}\oplus U^{-n+2}).$ 
By the spin representation, $[V\cdot V]\in u(E)\otimes_\R \w^2 T_M$ acts on $\psi$ to obtain 
$[V\cdot V]\cdot\psi \in \End(E)\otimes_\R  (U^{-n}\oplus U^{-n+2}).$
Then we have a definition of curvature of generalized connection $\D^\A$:
\bgn{definition}\label{def:Fa(psi)}[curvature of generalized connections]
 $\F_\A(\psi)$ of a generalized connection $\D^\A$ is defined by
\bgn{equation}{\F_\A}(\psi):=F_A\cdot\psi+ d^A(V\cdot\psi)+\frac12[V\, \cdot V]\cdot\psi
\end{equation}
\end{definition} 
$\F(\psi)$ is a globally defined section of $\End(E)\otimes (U^{-n}\oplus U^{-n+2})$ which is called 
the curvature of $\D^\A$.
\bgn{remark}
The ordinary curvature $F_A$ of a connection $D^A$ is defined to be the composition $d^A\circ d^A.$
However $\F_\a(\psi)$ is different from the composition 
$d^{\D_\A}\circ d^{\D_\A}$ which is not a tensor but a differential operator.
In fact, we have
\bgn{align}
d^{\D_\A}\circ d^{\D_\A}=&(d+\A)\circ (d+\A)=(d^A+ V)\circ (d^A+V)\\
=&F_A+ d^A\circ V + V\circ d^A+ V\circ V
\end{align} 
Then for $s\in \Gam(E)$, we have 
$$
(V\circ V)\cdot\psi =\frac 12[V\, \cdot V]\cdot\psi
$$
However, for $f\in C^\infty(M)$, $s\in \Gam(E)$, we have 
$$
(d^{\D_\A}\circ d^{\D_\A})(fs\otimes\psi)=f(d^{\D_\A}\circ d^{\D_\A})(s\otimes\psi)
-2\lan df, Vs\ran_{\tt}\cdot\psi
$$
Thus $(d^{\D_\A}\circ d^{\D_\A})$ is not a tensor but a differential operator.
\end{remark}
$\End(E)\otimes(U^{-n}\oplus U^{-n+2})$ is decomposed as 
$$\End(E)\otimes(U^{-n}\oplus U^{-n+2})=(\End(E)\otimes U^{-n})\oplus (\End(E)\otimes U^{-n+2})$$
We denote by $\pi_{U^{-n}}$ the projection from $\End(E)\otimes(U^{-n}\oplus U^{-n+2})$ to the component $\End(E)\otimes U^{-n}$.
 The line bundle $U^{-n}$ becomes the trivial complex line bundle by using the basis $\psi\in U^{-n}$.
Then $\End(E)\otimes_\C U^{-n}$ is identified with $\End(E).$ 
Then it follows from (\ref{eq:End(E)}) that we have 
 \bgn{align}
\End(E)
\otimes U^{-n}
=&u(E)\oplus \Herm(E,h)
\end{align}
We define $\pi_{\Herm}$ to be the projection to the component $\Herm(E,h)$
and we denote by $\pi^{\Herm}_{U^{-n}}$ the composition $\pi^{\Herm}\circ\pi_{U^{-n}}$.
Then we define $\K_\A(\psi)$ by 
$$\K_\A(\psi):=\pi_{U^{-n}}^{\Herm}\F_\A(\psi)\in \Herm(E,h)$$
\bgn{definition}\label{Einstein-Hermitian condition}[Einstein-Hermitian condition]
A generalized Hermitian connection $\D^\A$ is {\it Einstein-Hermitian} if $\D^\A$
satisfies the following: 
$$
\K_\A(\psi)=\lam id_E, \quad \text{\rm for a real constant }\lam
$$
\end{definition}
\bgn{remark}
If $\D^\A$ is an ordinary Hermitian connection $D^A$ over a K\"ahler manifold with
a K\"ahler form $\ome$, then 
the Einstein-Hermitian condition in Definition \ref{Einstein-Hermitian condition} coincides with the ordinary Einstein-Hermitian condition for $\psi=e^{\frac{\sqrt{-1}}2\ome}.$
In fact, the ordinary connection $D^A$ is Einstein-Hermitian connection if its curvature $F_A$ satisfies
$$\sqrt{-1}\Lam_\ome F_A=-\lam \id_E$$ or equivalently 
$$
\frac{\sqrt{-1}n F_A\w\ome^{n-1}}{\ome^n}= -\lam \id_E,
$$
where $\lam$ is a real constant.
The projection $\pi_{U^{-n}}$ of $F_A$ is given by 
$$
\pi_{U^{-n}}F_A =\frac{\lan F_A\cdot\psi, \, \ol\psi\ran_s}{\lan \psi, \,\ol\psi\ran_s}\psi
$$
since $\psi=e^{\frac{\sqrt{-1}}2 \ome}$, we have 
$$
\frac{\lan F_A\cdot\psi, \, \ol\psi\ran_s}{\lan \psi, \,\ol\psi\ran_s}=\frac{F_A\w (\sqrt{-1}\ome)^{n-1}}{(n-1)!}\frac{n!}{(\sqrt{-1}\ome)^n}=-\sqrt{-1} n\frac{F_A\w\ome^{n-1}}{\ome^n}
$$
Since $F_A$ is in $u(E)$, it follows that 
$$
-\sqrt{-1} n\frac{F_A\w\ome^{n-1}}{\ome^n}\in \Herm (E, h).
$$
Thus under the identification $U^{-n}\cong \C\psi$, we have 
$$
\pi^{\Herm}_{U^{-n}}F_A\cdot\psi =-\sqrt{-1} n\frac{F_A\w\ome^{n-1}}{\ome^n}\psi
$$
Hence $\pi^{\Herm}_{U^{-n}}F_A\cdot\psi=\lam \id _E$ is equivalent to 
$\sqrt{-1}\Lam_\ome F_A=-\lam \id_E$. 
\end{remark}
The unitary gauge group $U(E, h)$ acts on $\F_\A(\psi)$ by the adjoint action. Thus the Einstein-Hermitian condition is invariant under the action of the unitary gauge group.
Further our Einstein-Hermitian condition behaves nicely 
for the action of $b$-fields.
A real $d$-closed $2$-form $b$ acts on $\psi$ by 
$e^b\cdot\psi$ which is also a $d$-closed, nondegenerate, pure, spinor.
\bgn{definition}[$b$-field action of generalized connections]
Let $\D^\A=d+\A$ be a generalized connection of vector bundle $E$ over $(M, \psi).$  
A $d$-closed $2$-form $b$ acts on a generalized connection $\D^\A$ by 
$e^b\D^\A e^{-b}=d+e^b\A e^{-b}$, 
Then $e^b\A e^{-b}$ is given by 
$$
e^b\A e^{-b}=A+e^b Ve^{-b}=A+ad_b V +V,
$$
where $A+ ad_b V\in \End(E)\otimes T^*_M$ and 
$d+A+ \ad_b V$ is an ordinary connection of $E$ and 
$V\in \End(E)\otimes T_M$.
Note that $ad_b V$ is given by 
$$
\ad_bV =\sum_i V_i \otimes[b, v_i],
$$
for $V=\sum_i V_i \otimes v_i$, where $V_i\in u(E)$ and $v_i\in T_M$.
\end {definition}

\bgn{proposition}\label{prop: b-field action}
Let $\Ad_{e^b}\D^\A=d+\Ad_{e^b}\A$ be a generalized connection which is given by 
the $b$-field action on $\D^\A$. 
Then the curvature $\F_{\Ad_{e^b}\A}(\psi)$ is given by 
$$
\F_{\Ad_{e^b}\A}(\psi)=e^b\F_\A(e^{-b}\psi)
$$
 \end{proposition}
 \bgn{proof}
 The $b$-field action is given by 
  $\Ad_{e^b}\D^\A=d+A+\ad_bV + V,$ where $d+A+\ad_bV$ is a connection and its curvature is given by 
 $$
 F_A+d^A(\ad_b V)+\frac12[\ad_b V, \, \ad_b V]
 $$
 Then we have 
 \bgn{align}
\F_{\Ad_{e^b}\A}(\psi)= & (F_A+d^A(\ad_b V)+\frac12[\ad_b V, \, \ad_b V])\cdot\psi\\
+&d^A(V\cdot\psi)+[\ad_bV \cdot V]\cdot\psi+\frac12 [V\cdot V]\cdot\psi\\
\end{align}
Applying $\Ad_{e^b}V=V+\ad_b V$ and $e^bF_A e^{-b}=F_A$ and 
$[e^b Ve^{-b}\,\cdot\, e^b Ve^{-b}]=e^b[V\, \cdot\, V]e^{-b}$,  we have
\bgn{align}
\F_{\Ad_{e^b}\A}(\psi)=&F_A\cdot\psi+ d^A(e^b Ve^{-b}\cdot\psi)+
\frac12[e^b Ve^{-b}\,\cdot\, e^b Ve^{-b}]\cdot\psi\\
=&e^b\(F_A+ d^A(V\cdot\psi)+\frac12[V\cdot V]\)e^{-b}\cdot\psi\\
=&e^b\F_\A(e^{-b}\psi)
\end{align}
 \end{proof}
 Thus our Einstein-Hermitian condition is equivalent under  the action of $b$-field.
 \bgn{proposition}\label{prop: invariance of b-field action}
 Let $\D^\A$ be a generalized Hermitian connection of $E$ over $(M, \psi)$. Then  $\D^\A$ is Einstein-Hermitian over $(M,\psi)$ if and only if
 $\Ad_{e^b}\D^\A$ is an Einstein-Hermitian generalized connection 
 of $E$ over $(M, e^{b}\psi)$.
 \end{proposition}
 \bgn{proof}
 We denote by $\pi_{U^{-n}_{e^b\psi}}$ the projection to the component 
 $U^{-n}_{e^b\psi}:=e^b\cdot U^{-n}_{\J_\psi}.$ 
 From Proposition \ref{prop: b-field action}, we have 
 \bgn{align}
 \pi^{\Herm}_{U^{-n}_{e^b\psi}}\F_{\Ad_{e^b\A}}(e^b \psi)=&\pi^{\Herm}
 \frac{\lan \F_{\Ad_{e^b}}(e^b \psi), e^b\cdot\ol\psi\ran_s}{\lan e^b\psi, \,e^b\ol\psi\ran_s}\\
=&\pi^{\Herm}
 \frac{\lan e^b\F_{\A}(\psi), e^b\cdot\ol\psi\ran_s}{\lan e^b\psi, \,e^b\ol\psi\ran_s}
 \end{align}
 Since $\lan e^b\psi, \,e^b\ol\psi\ran_s=\lan \psi, \, \ol\psi\ran_s$, we have 
 $$
 \pi^{\Herm}_{U^{-n}_{e^b\psi}}\F_{\Ad_{e^b\A}}(e^b \psi)=
 \pi^{\Herm}
 \frac{\lan \F_{\A}(\psi), \,\cdot\ol\psi\ran_s}{\lan \psi, \,\ol\psi\ran_s}
=\pi^{\Herm}_{U^{-n}_{\J_\psi}}\F_\A(\psi).
 $$
 Thus $\pi^{\Herm}_{U^{-n}_{e^b\psi}}\F_{\Ad_{e^b\A}}(e^b \psi)=\lam\id_E$
 if and only if $\pi^{\Herm}_{U^{-n}_{\J_\psi}}\F_\A(\psi)=\id_E$. 
 Thus the result follows.
  \end{proof}
\subsection{The first Chern class of $E$ and $\tr\F_\A(\psi)$}
\bgn{theorem}\label{th: Chern class}
$\frac{-1}{2\pi\sqrt{-1}}\tr\F_\A(\psi)$ is a $d$-closed differential form on $M$ which is a representative of the class $[c_1(\det E)]\cup [\psi]\in H^\bullet(M)$.
\end{theorem}
\bgn{proof}As in Definition \ref{def:Fa(psi)}, $\tr\F_\A(\psi)$ is given by 
$$
\tr\F_\A(\psi):=\tr F_A\cdot\psi+ \tr d^A(V\cdot\psi)+\tr\frac12[V\, \cdot V]\cdot\psi
$$
We have  $\tr[V\cdot V]\cdot\psi=0$. 
By using local trivializations of $E,$ we have 
$$
\tr d^A(V\cdot\psi)=\tr \(d(V_\a\cdot\psi)+[A_\a\cdot V_\a]\)\cdot\psi 
=d (\tr V)\cdot\psi.
$$
Since $\tr V\cdot \psi$ is a globally defined form on $M,$ we have 
$$
\frac{-1}{2\pi\sqrt{-1}}[\tr\F_\A(\psi)]=\frac{-1}{2\pi\sqrt{-1}}[\tr F_A\cdot\psi]=[c_1(\det (E)]\cup[\psi].
$$
\end{proof}
\bgn{proposition}
Let $\D^\A$ be a generalized Einstein-Hermitian connection which satisfies $\pi^{\Herm}_{U^{-n}}\F_\A(\psi)=\lam \id_E$. 
Then $\lam$ is given in terms of the first Chern class $c_1(E)$ of $E$ 
and the class $[\phi]$ by 
$$
\frac{-\lam r}{2\pi\sqrt{-1}}\int_M\lan \psi, \, \ol\psi\ran_s =
\int_M \lan c_1(E)\w\psi ,\,\ol\psi\ran_s
$$
\end{proposition}
\bgn{proof}
It follows from Theorem \ref{th: Chern class} that we have 
$$
\frac{-1}{2\pi\sqrt{-1}}\int_M\tr\lan \F_\A(\psi), \, \ol\psi\ran_s =
\int_M\lan c_1(E)\w\psi, \, \ol\psi\ran_s.
$$
The $U^{-n}$-component of $\F_\A(\psi)$ is written as 
$$
\F_\A(\psi)=\lam \id_E\psi+ \xi\psi,
$$
where $\lam\in \R$ and $\xi\in u(E)$. 
Then we have 
\bgn{align*}
\frac{-1}{2\pi\sqrt{-1}}\tr\lan \F_\A(\psi), \, \ol\psi\ran_s
=&\frac{-1}{2\pi\sqrt{-1}}\tr \lan \lam \psi, \, \ol\psi\ran_s+ 
\frac{-1}{2\pi\sqrt{-1}}\tr\lan \xi \psi, \,\ol\psi\ran_s\\
=&\frac{-\lam r}{2\pi\sqrt{-1}}\lan \psi, \, \ol\psi\ran_s
+\frac{-1}{2\pi\sqrt{-1}}(\tr\xi)\lan \psi, \, \ol\psi\ran_s
\end{align*}
We have 
$$
\frac{\lan c_1(E)\w\psi, \, \ol\psi\ran_s}{\lan \psi, \, \ol\psi\ran_s}
=\frac{c_1(E)\w(\sqrt{-1}\ome)^{n-1}}{(n-1)!}\frac{n!}{(\sqrt{-1}\ome)^n}
=\frac{c_1(E)\w \ome^{n-1}}{\ome^n}\frac{n}{\sqrt{-1}}\in \sqrt{-1}\R
$$
Since $\xi \in u(E)$, $\tr\xi$ is pure imaginary. 
Thus we have the result.
\end{proof}
\section{Generalized holomorphic vector bundles}
Let $E$ be a complex vector bundle over a \complex manifold $(M, \J)$. 
{\it A generalized holomorphic structure} of $E$ is a differential operator 
$$
\ol\pa_\J : E\to E\otimes \ol \L_\J,
$$
which satisfies 
$$
\ol\pa_\J (f s) = s\otimes(\ol\pa_\J f)   + f(\ol\pa_\J s), \qquad \text{\rm for   }f\in C^\infty_M, s\in E.
$$
A complex vector bundle with a generalized holomorphic structure is called as
{\it a generalized holomorphic vector bundle}
\cite{Gua_2010}.
Then a generalized holomorphic structure $\ol\pa_\J$ is extended to be  
an operator $\ol\pa_\J : E \otimes \w^i\ol \L_\J \to E\otimes \w^{i+1}\ol\L_\J$ by 
$$
\ol\pa_\J (s \a) =(\ol\pa_\J s )\a+ s (\ol\pa_\J s),\qquad \text{\rm for   }
s\in E, \a\in \w^i\ol\L_\J.
$$
Then we obtain an elliptic complex which is the Lie algebroid complex.
$$
0\to E\to E\otimes \ol\L_\J \to  E\otimes \w^2\ol\L_\J\to \cdots \to  E\otimes \w^i\ol\L_\J \to 
E\otimes \w^{i+1}\ol\L_\J \to \cdots 0
$$
We denote by 
$C^{0,\bullet}=\{C^{0,i}\}$ the Lie algebroid complex  
$(E\otimes \w^i\ol{\L_\J}, \ol\pa_\J)$. 
The complex $C^{0,\bullet}$ is the deformation complex of generalized holomorphic structures on $E$. 
\subsection{The canonical generalized connection}
Let $(E, h)$ be a Hermitian vector bundle over a generalized complex manifold $(M,\J)$. 
We assume that $E$ admits a generalized holomorphic structure $\ol\pa_\J$.
Then there is a unique generalized Hermitian connection $\D^\A$ such that 
$$
\D^\A=\D^{1,0}+\ol\pa_\J,
$$
 where $\D^{1,0}: E\to E\otimes\L_\J$ denotes the $(1,0)$-component of $\D^\A$ with respect to $\J$. 
 We call $\D^\A$ as the canonical generalized connection of Hermitian vector bundle over a
 \complex manifold. 
 This is an analogue of the canonical connection of Hermitian vector bundle over a
 complex manifold. 
 In fact, $\D^\A$ is determined by 
 $$
 \pa_\J h(s_1, s_2) =h(\D^{1,0}s_1, s_2)+h(s_1, \ol\pa_\J s_2), \qquad s_1, s_2\in E.
 $$
\section{Curvature of line bundles over  generalized K\"ahler manifolds}
\subsection{The canonical generalized connection of the canonical line bundle over generalized K\"ahler manifolds}
Let $\J$ be a generalized complex structure on a manifold $M$ which gives the decomposition 
$$\w^\bullet T^*_M =\oplus_{i=0}^{2n} U^{-n+i}.$$
Then the canonical line bundle $K_\J$ of $\J$ is defined as $U^{-n}.$
A section of $K_\J=U^{-n}$ is a smooth differential form which is a nondegenerate, pure spinor. 
By using the action of $\ol \L_\J$ on $U_\J^{-n}$  given by spin representation , we have an identification  $U^{-n+1}_\J\cong K_\J\otimes \ol \L_\J.$
Thus the operator $\ol\pa_\J : U^{-n}_\J \to U^{-n+1}_\J$  is regarded as an operator $K_\J\to K_\J\otimes \ol\L_\J$ which satisfies 
$\ol\pa_\J(f s) =(\ol\pa_\J f)s+ f\ol\pa_\J s$. 
Since $\ol\pa_\J\circ \ol\pa_\J =0$, $\ol\pa_\J$ is a generalized holomorphic structure on $K_\J$. 
The canonical line bundle $K_\J$ of $\J$ admits trivializations $\{\phi_\a\}$ and transition functions $e^{\k_{\a\b}}$ relative to a cover $\{U_\a\}$, where $\phi_\a$ is a nondegenerate, pure spinor on $U_\a$ which induces the almost generalized complex structure $\J.$ Then the exterior derivative $d\phi_\a$ is given by 
$$
d\phi_\a= \eta_\a\cdot\phi_\a,
$$
where $\eta_\a$ denotes a real section of $\TT$ . Since $\eta_\a$  are real, they are uniquely determined. 
Since $d=\pa_\J+\ol\pa_\J$, the generalized holomorphic structure on $K_\J$ as before is given by connection  forms $\{\eta^{0,1}_\a\}, $ where 
$\eta_\a$ is decomposed into $\eta_\a^{1,0}+\eta_\a^{0,1}$ for 
$\eta_\a^{1,0}\in \L_\J$ and $\eta_\a^{0,1}\in \ol \L_\J$. 
Since $\eta_\a$ is  real, we have $\eta_\a^{1,0}=\ol {\eta_\a^{0,1}}$.
Let $\psi=e^{b+\sqrt{-1}}$  be a 
generalized complex structure given by a real symplectic form $\ome$ and a real $d$-closed $2$-from $b$.
We assume that the pair  $(\J, \psi)$ is a generalized K\"ahler manifold. 
A real function $\rho_\a$ is defined by 
$$
\lan \phi_a,\, \ol\phi_\a\ran_s =\rho_\a\lan \psi, \,\ol\psi\ran_s.
$$
Then we define $\A_\a$ by 
$$
\A_\a=\sqrt{-1}(-\J\eta_\a+\frac12\J d\log\rho_\a)\in \TT
$$
Then it is shown in \cite{Goto_2016} that 
$$
\A_\a\cdot\psi=\A_\b \cdot\psi+\sqrt{-1}d(\Im \k_{\a\b})\cdot\psi
$$
Since $e^{\sqrt{-1}\Im\k_{\a\b}}$ are regarded as transition functions of $K_\J$, 
we have 
\bgn{definition}\label{def: DA of KJ}
Thus $\A_\a$ gives a generalized connection $\D_\A$ on the canonical line bundle $K_\J$.
Then the curvature $\F_\A$ of $\D_\A$ is given by 
$$
\F_\A(\psi)=d(\A_\a\cdot\psi)=\sqrt{-1}d(-\J\eta_\a+\frac12\J d\log\rho_\a)\cdot\psi
$$
\end{definition}
Then $\F_\A(\psi)=(P+\sqrt{-1}Q)\cdot\psi$, where $-P$ coincides with the generalized Ricci curvature of 
a generalized K\"ahler manifold $(\J,\psi)$ in the paper" scalar curvature as moment map" \cite{Goto_2016}.
\bgn{proposition}\label{prop: generalized Chern connection of KJ}
Let $K_\J$ be the canonical line bundle over a generalized K\"ahler manifold 
$(M, \J, \J_\psi)$. 
We denote by $h$ the Hermitian metric of $K_\J$ which is given by 
$$
\lan \phi_\a, \, \ol\phi_\a\ran_ s= \rho_\a\lan \psi,\, \ol\psi\ran_s
$$
Then the generalized connection $\D^\A$ as in Definition \ref{def: DA of KJ}
 is the canonical generalized connection of the Hermitian generalized line bundle 
 $(K_\J, h)$.
\end{proposition}
\bgn{proof}
It suffices to show that the generalized connection $\D^\A$ preserves the Hermitian metric $h$ and $(0,1)$-component of $\D^{\A}$ coincides with the generalized holomorphic structure 
$\ol\pa_\J$ of $K_\J$. 
In fact, if we take $\phi_\a$ as a unitary basis of $K_\J$ with respect to $h$, then $\rho_\a=1$. Then we have 
$$
\A_\a=-\sqrt{-1}\J \eta_\a
$$
Since we take $\eta_\a$ are real, then $\sqrt{-1}\eta_\a$ are pure imaginary. The it follows that $\D^\A$ is a Hermitian generalized connection.
We also have 
$$
\ol\pa_\J\phi_\a=\eta^{0,1}\cdot\phi_\a
$$
The $(0,1)$-component of $\D^\A$ is given by 
$\eta^{0,1}$ with respect to local basis $\{\phi_\a\}$. 
Thus  $\ol\pa_\J $ is the $(0,1)$-component of $\D^\A$. 
Then the result follows.
\end{proof}
\bgn{proposition}
The curvature of the canonical generalized connection $\D^\A$ of the generalized Hermitian line bundle $(K_\J, h)$ over $(M ,\J, \J_\psi)$ 
is the moment map as in \cite{Goto_2016}.
\end{proposition}
\bgn{proof}
The result follows from Proposition \ref{prop: generalized Chern connection of KJ}
\end{proof}
\subsection{Einstein-Hermitian condition on generalized connections of line bundles and K\"ahler Ricci-solitons over ordinary K\"ahler manifolds}
Let $L$ be a Hermitian line bundle with a Hermitian metric $h$ over a K\"ahler manifold $(M, J ,\ome)$. 
We denote by $\J_J$ the \complex structure induced from the ordinary complex structure $J$.
Let $\D^\A$ be a generalized Hermitian connection $\D^\A$ which is written as $\D^\A:=d+\A=d+A+\sqrt{-1}v$, 
where $D^A=d+A$ is a Hermitian connection of $(L, h)$ and $v$ is a real vector filed. We have the decomposition
$\D_\A=\D^{1,0}_\A+\ol\pa_\A$ with respect to $\J_J$ and then the condition 
$\ol\pa_\A\circ\ol\pa_\A=0$ is equivalent to both equations
$$\ol\pa_A\circ\ol\pa_A=0 \,\text{\rm  and   }\ol\pa_A v^{1,0}=0.$$
Thus $\ol\pa_A$ is the ordinary  holomorphic structure on $L$ and
$v=v^{1,0}+v^{0,1}$ is a real holomorphic vector field with respect to $J$. 
Let $\psi$ be a $d$-closed, nondegenerate, pure spinor which is defined by 
$$
\psi =e^{c\ome+\sqrt{-1}\ome}.
$$
Note that we take a $b$-field to be $c\ome$, where $c$ is a real constant. 
Then $(\J_J, \psi)$ gives a generalized K\"ahler structure $(\J_J, J_\psi)$ as before. 
\bgn{proposition}\label{prop:line bundle}
Let $\D^\A:=D^A+\sqrt{-1}v$ be a generalized Hermitian connection of a Hermitian line bundle $(L, h)$ over generalized K\"ahler manifolds $(M, \J_J, \J_\psi)$. 
A Generalized Hermitian connection $\D^\A$ is Einstein-Hermitian if and only if 
 $(D^A, v)$ satisfies the following
$$
\Lam_\ome (F_A +c\sqrt{-1}L_v \ome) ={\lam} \sqrt{-1} ,
$$
where $\Lam_\ome$ denotes the contraction by K\"ahler form $\ome$.
\end{proposition}

\bgn{proof}
Let $\F_\A(\psi)$ be the curvature form of $\D^\A$ as in Definition \ref{def:Fa(psi)}.
Since $L$ is a line bundle, we have 
$[v\cdot v]=0$ and $[A\cdot v]=0$
we have 
\bgn{align}
\F_\A(\psi)
=&F_A\cdot\psi+\sqrt{-1}L_v\psi =F_A\w \psi+(c\sqrt{-1}L_v\ome-L_v\ome)\w\psi
\end{align}
The projection $\pi_{U_\psi^{-n}}$ is given by 
$$
\pi_{U_\psi^{-n}}\F_\A(\psi)=\frac{\lan \F_\A(\psi),\, \ol\psi\ran_s }{\lan \psi, \,\ol\psi\ran_s} \psi
$$
Since $\lan e^{-\ome}\F_\A(\psi), \, e^{-\ome}\ol\psi\ran_s =\lan \F_\A(\psi), \, \ol\psi\ran_s$, 
Thus we have 
\bgn{align}
\pi_{U_\psi^{-n}}\F_\A(\psi)
=&\frac1{\sqrt{-1}}\Lam_\ome (F_A+c\sqrt{-1}L_v\ome-L_v\ome)\psi.
\end{align}
Since $\Herm (L,h)$ is identified with the real part, the projection to $\Herm(L,h)$ is given  by 
$$
\pi^{\Herm}_{U^{-n}}\F_\A(\psi)=\frac1{\sqrt{-1}}\Lam_\ome (F_A+c\sqrt{-1}L_v\ome)
$$
Hence the Einstein-Hermitian condition $\pi^{\Herm}_{U^{-n}}\F_\A(\psi)=\lam \id$ is given by 
$$
\Lam_\ome (F_A+c\sqrt{-1}L_v\ome)=\lam\sqrt{-1}
$$
\end{proof}
\bgn{definition}[K\"ahler-Ricci soliton]
Let $D^A$ be the ordinary canonical connection of the anticanonical line bundle $(K^{-1}, h)$ over a K\"ahler manifold $(M, J, \ome)$ and $v$ a real holomorphic vector field on $(M, J)$,
where we denote by $h$ the canonical Hermitian metric on $K^{-1}$.
A pair $(D^A, v)$ is a K\"ahler Ricci soliton if $(D^A, v)$ satisfies the following 
$$
F_A+c\sqrt{-1}L_v\ome =\sqrt{-1}\ome, 
$$
where $F_A$ is given by the Ricci form $\text{\rm Ricc}_\ome$.
\end{definition}
Since $v$ is a real holomorphic and $\ol\pa^A\circ \ol\pa ^A=0$, $L$ is a generalized holomorphic vector bundle. Then we have the canonical generalized connection $\D^\A=D^A+\sqrt{-1}v$ of $L$.
Then a K\"ahler-Ricci soliton gives an Einstein-Hermitian generalized connection.
\bgn{proposition}
If $(D^A, v)$ is a K\"ahler-Ricci soliton, then the corresponding canonical generalized connection  $\D^\A=D^A+\sqrt{-1}v$ is Einstein-Hermitian over $(M, \J_J, \J_\psi)$. 
\end{proposition}
\bgn{proof}
The result follows from Proposition \ref{prop:line bundle}.
\end{proof}
\section{Curvature of generalized connections as moment map over a generalized K\"ahler manifold}
Let $(M^{2n}, \J, \J_\psi)$ be a compact generalized K\"ahler manifold of dimension $2n$, where $\J_\psi$ is a generalized K\"ahler structure induced from $\psi=e^{b+\sqrt{-1}\ome}$ for a real symplectic structure $\ome.$
The volume form $\vol_M$ on $M$ is defined by 
$$
\vol_M=i^{-n}\lan \psi,\, \ol\psi\ran_s.
$$
Let $(E, h)$ be a Hermitian vector bundle over $M$ with a Hermitian metric $h.$
Then we  define the followings: 
\bgn{align}
\GCon(E, h):=&\{\D^\A: \text{\rm generalized Hermitian connection of } (E,h)\}\\
\Con(E, h):=&\{ D^A:\text{\rm Hermitian connection of }(E,h)\}\\
U(E, h) :=& \text{\rm unitary gauge group}
\end{align}
The Lie algebra of the unitary gauge group $U(E,h)$ is denoted by $u(E)$ which is given by 
$$
u(E):=\{ \xi\in\End(E)\,|\, h(\xi s, s')=-h(s, \xi s'), \text{\rm for }s,s'\in \Gam(E)\, \}\\
$$
$\End(E)$ is decomposed as 
$$
\End(E)=u(E)\oplus \Herm(E,h),
$$
where 
$$
\Herm(E, h):=\{ \xi\in\End(E)\,|\, h(\xi s, s')=h(s, \xi s'), \text{\rm for }s,s'\in \Gam(E)\, \}.
$$
The space of Hermitian connections $\Con(E, h)$ is an affine space whose model is the vector space of sections $\Gam(u(E)\otimes T^*_M).$
The space of generalized Hermitian connections
$\GCon(E, h)$ is regarded as the cotangent bundle of $\Con(E, h)$ and 
$\GCon(E, h)$ is also an affine space whose model is the vector space of sections $\Gam(u(E)\otimes (\TT)).$
The gauge group $U(E, h)$ acts on both $\GCon(E, h)$ and $\Con(E, h)$.
The tangent space $T_\A\GCon(E, h)$ at $\D^{\A}$ is $\Gam(u(E)\otimes(\TT))$ and then each tangent is denoted by 
$\dot\A\in \Gam(u(E)\otimes(\TT)).$
The generalized complex structure $\J$ acts on $\Gam(u(E)\otimes(\TT))$ by 
$J(\xi\otimes e)=\xi \otimes \J e,$ for $\xi \in u(E)$ and $e\in \TT.$
Then
a complex structure $J_{\GCon}$ on the tangent space $T_{\A} \GCon (E,h)$ is given by the action of $\J$
$$J_{\GCon} : \dot\A\mapsto \J\dot\A.$$
Then it turns out that $J_{\GCon}$ is a complex structure on the affine space $\GCon(E, h).$
We define a symmetric point-wise coupling $\tr\lan\, ,\, \ran_{\tt}$ of $\Gam(u(E)\otimes (\TT))$ by 
$$
\tr\lan \xi_1\otimes e_1, \, \xi_2\otimes e_2\ran_{\tt}:=\tr(\xi_1\xi_2)\lan e_1, \, e_2\ran_{\tt} \in C^\infty(M)
$$
Then a global coupling $\tr\lan\, ,\,\ran_M$ is given by the integration over $M$
$$
\lan \xi_1\otimes e_1, \, \xi_2\otimes e_2\ran_M :=-\int_M\tr\lan \xi_1\otimes e_1, \, \xi_2\otimes e_2\ran_s\vol_M
$$
The another \complex structure $\J_\psi$ also acts on $\Gam(u(E)\otimes (\TT))$ by 
$\dot \A\mapsto \J_\psi\dot \A.$
Then a symplectic structure $\ome_{\GCon}$ on $\GCon(E,h)$ is defined by 
\bgn{equation}\label{sym}
\ome_{\GCon}(\dot{\A}_1, \dot\A_2)=-\int_M \tr \lan \J_\psi\dot\A_1, \dot\A_2\ran_{\tt} \vol_M,
\end{equation}
where $\dot{\A}_1, \dot{\A}_2\in \Gam(u(E\otimes (\TT)).$
We also have a Riemannian metric $g_{\GCon}$ on $\GCon$ by 
$$
g_{\GCon}(\dot\A_1, \dot\A_2)=-\int_M \tr\lan \h G\dot\A_1, \dot\A_2\ran_{\tt}\vol_M,
$$
where $\h G =-\J\circ\J_\psi$. 
Then it turns out that $(\GCon(E, h), g_{\GCon}, J_{\GCon}, \ome_{\GCon})$ is a K\"ahler manifold.
(Note that this is the ordinary K\"ahler manifold, yet $\GCon(E,h)$ is an infinite dimensional manifold.
We can introduce a Sobolev norm to make $\GCon(E,h)$ a Banach or Hilbert manifold. 
Since this is the standard method in the gauge theory and
we do not mention about the completion of $\GCon(E,h)$ in this paper.) 
In fact, since $\GCon(E, h)$ is an affine space, there exists a flat-torsion free connection and $J_{\GCon}$ and $\ome_{\GCon}$ and $g_{\GCon}$ are parallel. Thus $\ome_{\GCon}$ is closed and 
it implies that $(J_{\GCon}, \ome_{\GCon})$ is a K\"ahler structure on the infinite dimensional affine space $\GCon(E, h).$
Then the gauge group $\G(E, h)$ acts on $\GCon(E, h)$ preserving 
the K\"ahler structure $(J_{\GCon}, \ome_{\GCon})$.
Let $\D^\A$ be a generalized connection of $E$ and  
 $\F_\A(\psi)$ the curvature defined as in Definition \ref{def:Fa(psi)}. 
 We also define ${\F_\A}(\ol\psi)$ by
$$ {\F_\A}(\ol\psi):=F_A\cdot\ol\psi+ d^A(V\cdot\ol\psi)+\frac12[V\, \cdot V]\cdot\ol\psi
$$
Then we have
\begin{proposition}\label{prop: moment map}
Let  $\mu:\GCon(E, h)\to u(E)^*$ be the moment map on $\GCon(E, h)$ for the action of 
$U(E, h)$. Then the moment map $\mu$ is given by
$$
\lan\mu(\A),\xi\ran=\int_M \Im \,\,i^{-n}\tr \lan\xi\psi,\,\, \F_\A(\ol\psi)\ran_s,
$$
where $\mu(\A)\in u(E)^*$ denotes the value of $\mu$ at a generalized connection $\D^\A\in \GCon(E,h)$ and the coupling $\mu(\A)$ and
$\xi \in u(E)$ is denoted by $\lan\mu(\A),\xi\ran$ and $\Im$ is the imaginary part of $i^{-n}\tr \lan\xi\psi,\,\, \F_\A(\ol\psi)\ran_s.$
\end{proposition}
We need several Lemmas and Propositions for proof of Proposition \ref{prop: moment map}. 
Our proof of Proposition \ref{prop: moment map} is given after Lemma \ref{lem:keylemma3}.
\bgn{lemma}\label{lem:keylemma1}
$$
d\,\tr \lan\xi\cdot\psi,\,a\cdot\ol\psi\ran_s
=\tr\lan \D_\A\xi, \,\, a\cdot\ol\psi\ran_s-
\tr\lan \xi\cdot\psi,\,\, d^{\D_\A} (a\cdot\ol\psi)\ran_s
$$
where $\xi \in u(E)$ and $a\in u(E)\otimes(\TT)$
\end{lemma}
\bgn{proof}
We will calculate both sides by using local trivialization of $E$ over  $U_\a$ as in Section 3.
It follows from $d\sig=-\sig d$ that we have 
\bgn{align}
d\, \tr (\xi\cdot\psi\w\sig(a\cdot\ol\psi)) 
=&\tr(d\xi\cdot\psi\w\sig(a\cdot\ol\psi))-\tr(\xi\cdot\psi\w\sig(da\cdot\ol\psi))
\end{align}
Thus we have 
\bgn{align}
d\, \tr \lan\xi\cdot\psi,\,a\cdot\ol\psi\ran_s
=&\tr\lan d\xi\cdot\psi,\,a\cdot\ol\psi\ran_s-\tr\lan\xi\cdot\psi, \,d(a\cdot\ol\psi)\ran_s
\end{align}
Since $\lan e_i\cdot\psi, \,\,e_j\cdot\ol\psi\ran_s=
-\lan\psi,\,\, e_i\cdot e_j\cdot\ol\psi\ran_s$, using a local basis $\{e_i\}$ of $\TT$, we have 
\bgn{align*}
\tr\lan [\A_\a, \xi]\cdot\psi,\,\,a\cdot\ol\psi\ran_s=&
\sum_{i,j}\tr\lan [\A_{\a,i},\xi]e_i\cdot\psi,\,\, a_je_j\cdot\ol\psi\ran_s=\sum_{i,j}\tr\(\[\A_{\a,i},\xi]a_j\)
\lan e_i\cdot\psi, \,\,e_j\ol\psi\ran_s\\
=&-\sum_{i,j}\tr([\A_{\a,i},\xi]a_j)\lan \psi,\,\,e_i\cdot e_j\cdot\ol\psi\ran_s\\
=&\sum_{i,j}\tr(\xi[\A_{\a,i},a_j])\lan \psi,\,\,e_i\cdot e_j\cdot\ol\psi\ran_s=\sum_{i,j}\tr\lan \xi\psi, [\A_{\a,i}, a_j]e_i\cdot e_j\cdot\ol\psi\ran_s\\
=&\tr\lan \xi\psi, \,\, [\A_\a\cdot a ]\cdot\ol\psi\ran_s
\end{align*}
Hence we obtain
\bgn{align}
\tr\lan [\A_\a, \xi]\cdot\psi,\,\,a\cdot\ol\psi\ran_s=
\tr\lan \xi\psi, \,\, [\A_\a\cdot a] \cdot\ol\psi\ran_s
\end{align}
Thus we have 
\bgn{align}
\tr d(\xi\cdot\psi\w\sig(a\cdot\ol\psi)) 
=&\tr(d\xi\cdot\psi\w\sig(a\cdot\ol\psi))-\tr(\xi\cdot\psi\w\sig(da\cdot\ol\psi))\\
+&\tr\lan [\A_\a, \xi]\cdot\psi,\,\,a\cdot\ol\psi\ran_s
-\tr\lan \xi\psi, \,\, [\A_\a\cdot a] \cdot\ol\psi\ran_s\\
=&\tr\lan (\D_\A\xi)\cdot\psi, \,\, a\cdot\ol\psi\ran_s-
\tr\lan \xi\cdot\psi,\,\, d^{\D_\A}( a\cdot\ol\psi)\ran_s
\end{align}
\end{proof}
\bgn{lemma}\label{lem:Re}
$$\lan e_i, \,e_j\ran_{\tt}\vol_M=\Re \,i^{-n}\lan e_i\cdot\psi, \, e_j\cdot\ol\psi\ran_s,\qquad \text{\rm for  } e_i, e_j\in \TT$$
where $\Re\, i^{-n}\lan e_i\cdot\psi, \, e_j\cdot\ol\psi\ran_s$ denotes the real part of $i^{-n}\lan e_i\cdot\psi, \, e_j\cdot\ol\psi\ran_s.$
\end{lemma}
\bgn{proof}
We have 
$\lan e_i\cdot\psi, \, e_j\cdot\ol\psi\ran_s=-\lan e_j\cdot e_i\cdot\psi, \, \ol\psi\ran_s.$
Then applying $e_i\cdot e_j+e_j\cdot e_i=-2\lan e_i, \, e_j\ran_{\tt},$ we have 
\bgn{align}
2\lan e_i, \,e_j\ran_{\tt}\vol_M=&-i^{-n}\lan(e_i\cdot e_j+e_j\cdot e_i)\cdot\psi, \,\ol\psi\ran_s\\
=&i^{-n}\lan e_i\cdot \psi, \, e_j\cdot\ol\psi\ran_s+ i^{-n}\lan e_j\cdot\psi, \, e_i\cdot\ol\psi\ran_s\\
\end{align}
Since $\lan \a, \, \b\ran_s=(-1)^n\lan \b,\a\ran_s $ for $\a, \b\in \w^\bullet T^*_M,$ we 
 have 
$$
\ol{i^{-n}\lan e_i\cdot\psi, \, e_j\cdot\ol\psi\ran_s }=(-1)^ni^{-n}\lan e_i\cdot\ol\psi,\,e_j\cdot\psi\ran_s
=i^{-n}\lan e_j\cdot\psi, e_i\cdot\ol\psi\ran_s
$$
Thus we have 
$$
2\lan e_i, \, e_j\ran_{\tt}\vol_M=i^{-n}\lan e_i\cdot\psi, \, e_j\cdot\ol\psi\ran_s+\ol{i^{-n}\lan e_i\cdot\psi, \, e_j\cdot\ol\psi\ran_s}
$$
Hence we have 
$$
\lan e_i, \,e_j\ran_{\tt}\vol_M=\Re \,i^{-n}\lan e_i\cdot\psi, \, e_j\cdot\ol\psi\ran_s.
$$
\end{proof}
We also have
\bgn{lemma}\label{lem:Im}
$$
\lan \J_\psi e_i, \, e_j\ran_{\tt}\vol_M =-\Im\, i^{-n}\lan e_i\cdot\psi, \, e_j\cdot\ol\psi\ran_s
$$
where $\Im\, i^{-n}\lan e_i\cdot\psi, \, e_j\cdot\ol\psi\ran_s$ denotes the imaginary part of 
$ i^{-n}\lan e_i\cdot\psi, \, e_j\cdot\ol\psi\ran_s.$
\end{lemma}
\bgn{proof}
Applying Lemma \ref {lem:Re}, we have 
\bgn{align}
\lan \J_\psi e_i, \, e_j\ran_{\tt}\vol_M=\Re \,i^{-n}\lan (\J_\psi e_i)\cdot\psi, \, e_j\cdot\ol\psi\ran_s
\end{align}
By using $\J_\psi$,  $e_i$ is decomposed into $e^{1,0}_\psi+ e^{0,1}_\psi$, where $e^{1,0}_\psi \in \L_{\J_\psi}$ and 
$e^{0,1}_\psi \in \ol{\L_{\J_\psi}}.$ 
Since $\J_\psi e^{0,1}_\psi =\sqrt{-1} e^{0,1}_\psi$ and $e^{1,0}_{\psi} \cdot\psi =0,$  we have 
\bgn{align}
(\J_\psi e_i)\cdot\psi=\sqrt{-1}e^{0,1}_\psi \cdot\psi =\sqrt{-1}e_i\cdot\psi
\end{align}
Thus we have 
$$
\lan \J_\psi e_i, \, e_j\ran_{\tt}=\Re \,\,i^{-n} \sqrt{-1}\lan e_i\cdot\psi, \, e_j\cdot\ol\psi\ran_s
=-\Im \,\,i^{-n}\lan e_i\cdot\psi, \, e_j\cdot\ol\psi\ran_s
$$
\end{proof}
\bgn{lemma}\label{lem:keylemma2}Let $\D^\A$ be a generalized Hermitian connection. 
Then the symplectic structure $\ome_{\GCon}$ at the tangent space $T_{\D^\A}\GCon(E,h)$ as before is given by 
$$
\ome_{\GCon}(\D_\A\xi, a)=
\int_M \Im\, i^{-n}\tr\lan \xi\psi,\,\, d^{\D_\A}( a\cdot\ol\psi)\ran_s,
$$
for $\xi \in u(E)$ and $a\in u(E)\otimes (\TT)$.
\end{lemma}
\bgn{proof}
As in (\ref{sym}),
the symplectic structure $\ome_{\GCon}$ on $\GCon(E,h)$ is defined by 
\bgn{equation}
\ome_{\GCon}(\dot{\A}_1, \dot\A_2)=-\int_M \tr \lan \J_\psi\dot\A_1, \dot\A_2\ran_{\tt} \vol_M,
\end{equation}
where $\dot{\A}_1, \dot{\A}_2\in \Gam(u(E\otimes (\TT)).$
By using a local real basis $\{e_i\}$ of $\TT$, we have  $\dot\A_1=\sum_i \dot\A_{1,i}e_i$ and $\dot\A_2=\sum_j\dot\A_{2,j}e_j.$
From Lemma \ref {lem:Im},  we have 
\bgn{align*}
\tr \lan \J_\psi\dot\A_1, \dot\A_2\ran_{\tt} \vol_M=&\sum_{i,j}\tr(\dot\A_{1,i}\dot\A_{2,j})
\lan\J_\psi e_i, \, e_j \ran_s\vol_M\\
=&-\sum_{i,j}\Im \, i^{-n}\tr(\dot\A_{1,i}\dot\A_{2,j})\lan e_i\cdot\psi, \, e_j\cdot\ol\psi\ran_s\\
=&-\Im \,i^{-n}\tr\lan\A_1\cdot\psi, \, \A_2\cdot\ol\psi\ran_s
\end{align*}
Then
the symplectic structure $\ome_{\GCon}$ is given by 
$$
\ome_{\GCon}(\dot\A_1, \dot\A_2)=\int_M \Im \, i^{-n} \tr \lan \dot\A_1\cdot\psi, \,\,\dot\A_2\cdot\ol\psi\ran_s
$$
for $\dot \A_1, \dot\A_2\in T_{\D^\A}\GCon(E,h)$.
Substituting $\dot\A_1=\D_\A\xi$ and $ \dot\A_2=a$, we have  
\bgn{align}
\ome_{\GCon}(\D_\A \xi, a)=&\int_M\Im\, i^{-n}\tr\lan \D_\A\xi\cdot\psi,\,\, a\cdot\ol\psi\ran_s
\end{align}
From Lemma \ref{lem:keylemma1} and Stokes' theorem, we have 
\bgn{align}
\ome_{\GCon}(\D_\A \xi, a)=&\int_M \Im\, i^{-n} \tr\lan \xi\psi,\,\, d^{\D_\A}( a\cdot\ol\psi)\ran_s
\end{align}
\end{proof}
Let 
$\D^{\A_t}=D^{A_t}+ V_t$ be a $1$-parameter family of  generalized Hermitian connections which consists of Hermitian connections $D^{A_t}$ and sections $V_t\in u(E)\otimes T_M. $
The infinitesimal deformations of $\D^{\A_t}$ at $\D^\A:=\D^{\A_0}$ is denoted by  $\dot{\A}:=\frac d{dt}\D^{\A_t}|_{t=0}\in u(E)\otimes (\TT),$ where $\dot{\A}=\dot A+\dot V$ and 
$\dot A\in u(E)\otimes T^*_M$ and $\dot V\in u(E)\otimes T_M.$ 
The operator $d^{\A}$ acts on $\dot A\cdot\psi\in u(E)\otimes U^{-n+1}$ to obtain 
$d^\A(\dot A\cdot\psi)\in \End(E)\otimes (U^{-n}\oplus U^{-n+2})$.
We also have the differential $\frac d{dt}\F_{\A_t}(\psi)|_{t=0}.$
Then the difference between $d^{\A}(\dot \A\cdot\psi)$ and $\frac d{dt}\F_{\A_t}(\psi)|_{t=0}$ is given by 
\bgn{proposition}\label{key lem 1}
$$d^{\A}(\dot \A\cdot\psi) =\frac d{dt}\F_{\A_t}(\psi)|_{t=0}+[V\cdot\dot A]\cdot\psi-[\dot A\cdot V]\cdot\psi
$$
\end{proposition}
\begin{proof} We shall calculate both sides in terms of 
local trivializations of $E$ of $\{U_\a\}$ which  are given by a unitary basis as before. 
The ordinary covariant derivative $d^A$ gives
$$
d^A(V\cdot\psi )=d_\a(V_\a\cdot\psi)+[A_\a\, \cdot V_\a]\cdot\psi,
$$
where $A_\a, V_\a\in u(E)|_{U_\a}$.
Hence we have 
\bgn{align}
\frac d{dt} d^{A_t}(V_t\cdot\psi )|_{t=0} =&d_\a (\dot V_\a\cdot\psi)+[\dot A_\a\cdot V_\a]\cdot\psi+
[A_\a \cdot \dot V_\a]\cdot\psi\\
=&d^A (\dot V\cdot\psi)+[\dot A\,\cdot V]\cdot\psi
\end{align}
It follows from $[\dot V\cdot V]=[V\cdot\dot V]$ that we obtain
$$
\frac d{dt}[V_t\, \cdot V_t]\cdot\psi =2[\dot V\, \cdot V]\cdot\psi
$$
Thus we have
\bgn{align}
\frac d{dt}\F_{\A_t}|_{t=0}=d_A(\dot A\cdot\psi)+d^A (\dot V\cdot\psi)+[\dot A\,\cdot V]\cdot\psi
+[\dot V \,\cdot V]\cdot\psi.
\end{align}
From (\ref{dDA(acdotpsi)}),
the right hand side $d^{\A}(\dot \A\cdot\psi)$ is given by 
\bgn{align}
d^{\A}(\dot{\A}\cdot\psi)=&d(\dot{\A}\cdot\psi)+[\A\, \cdot \dot{\A}]\cdot\psi\\
=&d(\dot{A}\cdot\psi)+d(\dot{V}\cdot\psi)\\
+&[A\, \cdot \dot{A}]\cdot\psi+
[A\, \cdot \dot{V}]\cdot\psi+[V\, \cdot \dot{A}]\cdot\psi+[V\, \cdot \dot{V}]\cdot\psi\\
=&d^A(\dot A\cdot\psi) +d^A(\dot V\cdot\psi)+[V\, \cdot \dot{A}]\cdot\psi+[V\, \cdot \dot{V}]\cdot\psi
\end{align}
Applying 
$[V\, \cdot \dot{V}]\cdot\psi=[\dot V\, \cdot {V}]\cdot\psi$ again, we obtain
$$
d^{\A}\dot \A =\frac d{dt}\F_{\A_t}(\psi)|_{t=0}+[V\cdot\dot A]\cdot\psi-[\dot A\cdot V]\cdot\psi
$$
\end{proof}
\bgn{proposition}\label{prop:key propsition}
Let $\F_\A(\ol\psi)$ be as in Definition \ref{def:Fa(psi)}
$$
{\F_\A}(\ol\psi):=F_A\cdot\ol\psi+ d^A(V\cdot\ol\psi)+\frac12[V\, \cdot V]\cdot\ol\psi
$$
Let $\D^{\A_t}$ be one parameter family of generalized Hermitian connection and $\dot \A$ the derivative of $\D^{\A_t}$ at $\D^\A=\D^{\A_0}$, that is,  
$\frac d{dt}\dot\D^{\A_t}|_{t=0}.$
Then we have the following
$$
\Im\, i^{-n}\tr \lan \xi\cdot\psi, \,(d^\A \dot\A)\cdot\ol\psi\ran_s=
\Im\, i^{-n}\tr \lan \xi\cdot\psi, \,\frac d{dt}\F_{\A_t}(\ol\psi)|_{t=0}\ran_s
$$
for every $\xi \in u(E)$, where $\frac d{dt}\F_\A(\ol\psi)|_{t=0}$ denotes the derivative at $t=0$.
\end{proposition}
\bgn{proof}
It follow from Lemma \ref{key lem 1} that we have 
\bgn{align}
&\Im\, i^n\tr\lan \xi\cdot\psi, \,d^\A( \dot\A\cdot\ol\psi)\ran_s-
\Im\, i^n\tr\lan \xi\cdot\psi, \,\frac d{dt}\F_\A(\ol\psi)|_{t=0}\ran_s\\
=&\Im \, i^n\tr\lan\xi\cdot\psi,\, [V\cdot\dot A]\cdot\ol\psi-[\dot A\cdot V]\cdot\ol\psi\ran_s
\end{align}
Then $V=\sum_i V_i e_i$ and $\dot\A=\sum \dot A_j e_j$, where $e_i, e_j\in \TT$ and $V_i, A_j\in u(E)$. 
Thus we have 
\bgn{align}
[V\cdot\dot A]\cdot\ol\psi-[\dot A\cdot V]\cdot\ol\psi=&\sum_{i,j}[V_i, \dot A_j]e_i\cdot e_j\cdot\ol\psi -
[\dot A_j, V_i]e_j\cdot e_i\cdot\ol\psi\\
=&\sum_{i,j}[V_i, \dot A_j]e_i\cdot e_j\cdot\ol\psi +[ V_i,\dot A_j]e_j\cdot e_i\cdot\ol\psi\\
=&\sum_{i,j}[V_i, \dot A_j](e_i\cdot e_j+e_j\cdot e_i)\cdot\ol\psi\\
=&-2\sum_{i,j}[V_i, \dot A_j]\lan e_i, \, e_j\ran_{\tt}\ol\psi
\end{align}
where $[V_i, \dot A_j]\in u(E)$ and $,\lan e_i, \, \, e_j\ran_{\tt}$ is a real function.
Thus it follows 
\bgn{align*}
\Im \, i^n\tr\lan\xi\cdot\psi,\, [V\cdot\dot A]\cdot\ol\psi-[\dot A\cdot V]\cdot\ol\psi\ran_s
=&-2\sum_{i,j}\Im \, i^n\tr\lan \xi \cdot\psi, \, [V_i, \dot A_j]\lan e_i, \, e_j\ran_{\tt}\ol\psi
\ran_s\\
=&-2\Im \(i^n\sum_{i,j}\tr(\xi  [V_i, \dot A_j])\lan e_i, \, e_j\ran_{\tt}\lan \psi, \, \ol\psi\ran_s \)
\end{align*}
Since $\xi,[V_i, \dot A_j]\in u(E)$, it follows 
$\tr(\xi  [V_i, \dot A_j])$ is real and $i^n\lan \psi, \, \ol\psi\ran_s$ is real also, thus
we have 
$$
\Im \, i^n\tr\lan\xi\cdot\psi,\, [V\cdot\dot A]\cdot\ol\psi-[\dot A\cdot V]\cdot\ol\psi\ran_s
=0
$$
Then the result follows. 
\end{proof}
Then we have 
\bgn{lemma}\label{lem:keylemma3}
Let $\mu(\A)\in u(E)^*$ be a value of the moment map $\mu$ at a generalized Hermitian connection $\D^\A$
and $a\in \Gam(M, u(E)\otimes (\TT))$ a tangent of $\GCon(E,h)$ at $\D^\A$. 
Then the moment map $\mu$ is given by 
$$
\lan d\mu(\A)a, \xi\ran=\int_M \Im\, i^n \tr \lan \xi\psi, \,\,   \frac d{dt}\F_{\A_t}(\ol\psi)|_{t=0}   \ran_S
$$
for $\xi \in u(E)$.
\end{lemma}
\bgn{proof}
We define a one-parameter family $\D^{\A_t}$ by $\D_\A + at$. 
Then $\frac d{dt}\mu(\D_{\A_t})|_{t=0}$ is denoted by $d\mu(\A)(a).$
From the definition of the moment map $\mu$,  we have 
$$
\lan d\mu(\A)(a), \xi \ran_s = \ome_{\GCon}(\D^\A\xi , \, a)
$$
From Lemma \ref
{lem:keylemma2}, we have 
$$
\lan d\mu(\A)(a), \xi \ran_s =\ome_{\GCon}(\D_\A\xi, a)=
\int_M \Im\, i^{-n}\tr\lan \xi\psi,\,\, d^{\D_\A}( a\cdot\ol\psi)\ran_s,
$$
It follows from Proposition \ref{prop:key propsition} that we have 
$$
\lan d\mu(\A)(a), \xi \ran_s=
\int_M\Im\, i^{-n}\tr \lan \xi\cdot\psi, \,\frac d{dt}\F_{\A_t}(\ol\psi)|_{t=0}\ran_s
$$
\end{proof}
{\noindent\sc Proof of Proposition \ref{prop: moment map}}. 
It follows from Lemmas \ref{lem:keylemma1}, Lemma \ref{lem:keylemma2} and Lemma \ref{lem:keylemma3}
that we have 
$$
\lan \mu (\A)(a), \xi \ran=\int_M \Im \,i^{-n}\lan \xi\cdot\psi,\, \F_\A(\ol\psi)\ran_s
$$
Thus we obtain the result. 
\qed \par\medskip
We consider the coupling between $u(E)$ and $\Herm(E,h)$ by 
$$
\int_M\Im \,i^{-n}\tr\lan \xi_1\cdot\psi, \, \sqrt{-1}\xi_2\cdot\ol\psi\ran_s
$$
for $\xi_1\in u(E)$ and $\sqrt{-1}\xi_2 \in \Herm (E,h)$. 
Since the coupling is perfect, it gives an identification $u(E)^*$ with $\Herm(E,h)$.
Then we have
\bgn{proposition}\label{prop:moment map2}
Under the identification $u(E)^*$ with $\Herm(E,h)$ as before, the moment map $\mu$ is given by 
$$
\mu(\A) =\K_\A(\psi):=\pi^{\Herm}_{U^{-n}}\F_\A(\psi)
$$
\end{proposition}
{\noindent\sc Proof of Proposition \ref{prop:moment map2}}.
Let $(\xi_1+\sqrt{-1}\xi_2)\cdot\psi$ be an element of $\End(E)\otimes U^{-n}$, where 
$\xi_1, \xi_2\in u(E)$.
Since $\xi \in u(E)$ and $i^{-n}\lan \psi, \, \ol\psi\ran_s$ is real and 
$\tr (\xi\xi_i)$ are also real for $i=1,2$,  we have 
\bgn{align}
\Im\, i^{-n}\tr\lan \xi \cdot\psi, \,\, (\xi_1+\sqrt{-1}\xi_2)\cdot\ol\psi\ran_s
=&\Im\, \tr(\xi\xi_1)i^{-n}\lan \psi, \, \ol\psi\ran_s+\Im\, \tr(\xi\xi_2)i^{-n}i\lan \psi, \, \ol\psi\ran_s\\
=&  \Im\, \tr(\xi\xi_2)i^{-n}i\lan \psi, \, \ol\psi\ran_s.\\
\end{align}
Since $\pi^{\Herm}(\xi_1+\sqrt{-1}\xi_2)\cdot\psi=\sqrt{-1}\xi_2\cdot\psi$, 
 we have 
$$
\Im\, i^{-n}\lan \xi\cdot\psi, (\xi_1+\sqrt{-1}\xi_2)\cdot\ol\psi\ran_s =
\Im\, i^{-n}\lan \xi \cdot\psi, \,\, \pi^{\Herm}(\xi_1+\sqrt{-1}\xi_2)\cdot\ol\psi\ran_s
$$
Hence we obtain 
$$
\Im\, i^{-n}\lan \xi\cdot\psi, \,\F_\A(\ol\psi)\ran_s=\Im \, i^{-n}\lan \xi \cdot\psi, \, 
\pi^{\Herm}_{U^{n}}\F_\A(\ol\psi)\ran_s
$$
Thus under the identification $u(E)^*\cong \Herm (E,h)$ by the coupling 
$\Im\, \tr\lan \, ,\,\ran_s$, the moment map $\mu$ is given by 
$$
\mu(\A)=\pi^{\Herm}_{U^{-n}}\F_\A(\psi)
$$
\qed
\bgn{proposition}\label{prop:Kahler quotient1}
Let $\lam$ be a real constant which is determined in terms of the $1$-st Chern class of $ E$
and the class $[\psi]$ as in the Einstein-Hermitian condition.
Then the quotient of the set of generalized connections satisfying the Einstein-Hermitian condition divided by the action of the unitary gauge group $U(E,h)$ is 
the K\"ahler quotient $\mu^{-1}(\lambda\id_E)/U(E, h)$. 
\end{proposition}
\bgn{proof}
The result immediately follows from Proposition \ref{prop: moment map} and Proposition \ref{prop:moment map2}.
\end{proof}
By  the decomposition $(\TT)^\C=\L_\J\oplus \ol\L_\J$ in terms of $\J,$ 
we  have the following decomposition 
$$E\otimes(\TT)=(E\otimes\L_\J )\oplus (E\otimes \ol\L_\J). $$
Then a generalized connection is also decomposed into 
$\D^\A=\D^{1,0}+\ol\pa^A_\J,$ where
$\D^{1,0}: \Gam(E)\to \Gam(E\otimes \L_\J)$ and $\ol\pa^\A_\J: \Gam(E)\to \Gam (E\otimes \ol \L_\J).$
Then $\ol \pa^\A_\J$ is extended to the operator $\Gam (E\otimes \w^p\ol\ L_\J)\to \Gam (E\otimes \w^{p+1}\ol \L_\J)$ for $p=0, 1,\cdots 2n$.
As in before, If $\ol\pa^\A_\J\circ \ol\pa^\A_\J=0$, then $(E, \ol\pa^\A_\J)$ is 
{\it generalized holomorphic vector bundle.}
\bgn{definition}
We define $\GCon(E, h)^{hol}$ to be a subset of $\GCon(E, h)$ consisting of $\D_\A$ such that 
$\ol\pa_\J^\A\circ \ol\pa_\J^\A=0$, that is, 
$$
\GCon(E, h)^{hol}:=\{\, \D^\A\in \GCon(E,h)\, |\, \ol\pa_\J^\A\circ \ol\pa_\J^\A=0\, \}
$$
\end{definition}
\bgn{remark}
This is an analogue of the notion of Hermitian connections associated with holomorphic vector bundles which are defined by 
$\ol\pa^A :\Gam(E) \to \Gam (E\otimes T^{0,1})$ satisfying $\ol\pa^A\circ \ol\pa^A=0$.
We also defined  $\Con(E, h)^{hol}$ by 
$$
\Con(E, h)^{hol}:=:=\{\, D^A\in \Con(E,h)\, |\, \ol\pa^A\circ \ol\pa^A=0\, \}
$$
\end{remark}
\bgn{proposition}\label{prop: subKahler mfd}
$\GCon(E, h)^{hol}$ is a complex submanifold of $\GCon(E, h)$ which inherits the K\"ahler structure.
\end{proposition}
\bgn{proof} The space $\Gam (M, u(E)\otimes_\R(\TT))$ is identified with $\Gam(M, \End(E)\otimes_\C 
\ol{\L_\J})$. Then
the tangent space of $\GCon(E,h)^{hol}$ at $\D^\A$ is the subspace of $\Gam (M, \End(E)\otimes\ol{\L_\J})$ which is given by 
$$
T_{\D^\A}\GCon(E, h)^{hol}:=\{ \dot \A \, |\, \ol\pa_\J^\A \dot \A^{0,1}=0\,\},
$$
where $\dot A=\dot \A^{1,0}+\dot \A^{0,1}$ and $\dot \A^{1,0}\in \Gam(M, \End(E)\otimes \L_\J)$ and $\dot \A^{0,1}\in \Gam(M, \End(E)\otimes \ol{\L_\J})$.
Thus $T_{\D^\A}\GCon^{hol}(E, h)$ is a complex subspace of $T_{\D^\A}\GCon(E, h)\cong \Gam(M,\End(E)\otimes \ol{\L_\J})$.
Since $J_{\GCon}\dot\A^{0,1}=\sqrt{-1}\dot\A^{0,1}$, $T_{D^A}{\GCon}^{hol}(E, h)$ is invariant under the action of $J_{\GCon}$. 
Then the result follows.
\end{proof}
The unitary gauge group $U(E,h)$ also preserves the K\"ahler structure on $\GCon(E, h)^{hol}.$
\bgn{theorem}\label{th:main theorem1}
The moment map of $\GCon(E, h)^{hol}$ is given by the restriction of the map $\mu$ as in Proposition \ref{prop: moment map}
 to the $\GCon(E, h)^{hol}$
\end{theorem}
\bgn{proof}
The result follows from Proposition \ref{prop: subKahler mfd}.
\end{proof}
We define the set of Einstein-Hermitian generalized connections $ {\GCon(E, h)}^{hol}_{EH}$  by 
$$
 {\GCon(E, h)}^{hol}_{EH}=
\{\, \D_\A\in \GCon(E, h)^{hol}\, |\, \K_\A(\psi)=\lam \id_E\, \}
$$
The moduli space of Einstein-Hermitian generalized connections ${{\cal GM}(E, h)}^{hol}_{EH}$ is defined to be the quotient 
of $ {\GCon(E, h)}^{hol}_{EH}$ by the action of the gauge group $U(E, h)$, that is, 
$$
{{\cal GM}(E, h)}^{hol}_{EH}= {\GCon(E, h)}^{hol}_{EH}/\G(E, h)
$$
\bgn{theorem}\label{th:main theorem 2}
The moduli space ${{\cal GM}(E, h)}^{hol}_{EH}$ arises as 
a K\"ahler quotient and the smooth part of the moduli space ${{\cal GM}(E, h)}^{hol}_{EH}$ inherits a K\"ahler structure.
\end{theorem}
\bgn{proof}
The result follows from Proposition \ref{prop:Kahler quotient1} and Proposition \ref{th:main theorem1}\footnote{We can construct a slice of ${\GCon(E, h)}^{hol}_{EH}$ for the action of $U(E,h)$ which becomes coordinates of the smooth part of the moduli space. This is the standard construction of Kuranishi space. In fact, as shown in next section, deformations of Einstein-Hermitian generalized connections are controlled   by an elliptic complex. }
\end{proof}
\section{Deformation complex of Einstein-Hermitian generalized connections}
Let $\D^\A$ be an Einstein-Hermitian generalized connection of a Hermitian vector bundle $(E, h)$ over a generalized K\"ahler manifold $(M, \J,\J_\psi)$.
We denote by $B^0$ the vector bundle $u(E)$ of 
skew-symmetric endomorphisms of $E$. 
Let $B^1$ be the vector bundle $u(E)\otimes_\R(\TT)$. 
The infinitesimal deformations of generalized Hermitian connections are given by global section of the bundle $B^1,$ that is, the tangent space of $\GC(E, h)$ at $\D_\A$
is the space of sections of $B^1.$
The unitay gauge group $U(E,h)$ acts on $\GCon(E, h)$ which gives 
infinitesimally the action of the Lie algebra $u(E):=B^0$ on the space $\GCon(E, h).$ The infinitesimal action is given by
the covariant derivative $d^{\D_\A} :B^0\to B^1$, that is, 
the tangent space of the action of  the gauge group is 
the image of $d^{\D_\A}(B^0).$
The complexified space $B^1_\C:=(u(E)\otimes(\TT))_\C$ is decomposed into 
$(\End(E)\otimes\L_\J)\oplus(\End(E)\otimes \ol\L_\J).$
Then an element $\xi\otimes e\in B^1_\C$ is given by 
$\xi\otimes e=\xi\otimes e^{1,0}+\xi\otimes e^{0,1}$, 
where $\xi\otimes e^{1,0}\in \End(E)\otimes \L_\J$ and $\xi\otimes e^{0,1}\in \End(E)\otimes \ol \L_\J.$
The action of $\tau$ on the space $B^1_\C$ is defined as 
$$
\tau (\xi\otimes e)=\xi^*\otimes\ol e
$$
where $\xi^*$ denotes the adjoint of $\xi\in\End(E)$ with respect to $h$. 
Then the $B^1$ is given by the real part of $B^1_\C$ which is
the skew-Hermitian part of $B^1_\C$, that is 
$B^1:=\{ \xi\otimes e \in B^1_\C\, |\, \tau(\xi\otimes e)=-\xi\otimes e\, \}.$
Then $\End(E)\otimes \ol \L_\J$ is identified with $B^1$ by 
$$
\xi\otimes e^{0,1}\mapsto -\xi^*\otimes \ol{e^{0,1}}+\xi\otimes e^{0,1}
$$
Under the identification $B^1$ with $\End(E)\otimes \ol{\L_\J}$,
the operator $\ol\pa_\J^\A :\End(E)\otimes \ol\L_\J\to \End(E)\otimes \w^2\ol \L_\J$  
is regarded as an operator 
form $B^1$ to $\End(E)\otimes \w^2\ol{\L_\J}.$ 
Each $\xi\otimes e\in \End(E)\otimes(\TT)$ acts on $\psi$ by 
$\xi\otimes e\cdot\psi \in \End(E)\otimes U^{-n+1}$. 
Then by using 
$$d^{\D_\A} : \End(E)\otimes U^{-n+1}\to \End(E)\otimes( U^{-n}\oplus U^{-n+2}),$$ we have $d^{\D_\A}(\xi\otimes e\cdot\psi)\in\End(E)\otimes( U^{-n}\oplus U^{-n+2}).$
The map $\K :B^1\to \Herm(E, h)$ is defined by 
$$
K(\xi\otimes e ): = \pi_{ U^{-n}}^{\Herm}d^\A ((\xi\otimes e)\cdot\psi) \in \Herm (E)
$$
where $\pi_{U^{-n}}^{\Herm}$ denotes the projection to $\Herm(E,h)$ as in Section 3.2.

We define $B^2_{\prim}$ to be the sum $\Herm(E, h)\oplus (\End(E)\otimes\w^2\ol\L_\J)$. 
Then we have the map $K\oplus \ol\pa_\J^\A: B^1\to B^2_{\prim}$ which is given by 
$$
K\oplus \ol\pa_\J^\A(\xi\otimes e):= K(\xi\otimes e)+\ol\pa_\J^\A(\xi\otimes e^{0,1})
$$
We also define $B^i$ be  $\End(E)\otimes\w^i\ol \L_\J$ for $i\geq 3$.
Then we have the following complex $B^\bullet$: 
$$
0\to B^0\overset {d^{\D_\A}}\longrightarrow B^1\overset{K\oplus \ol\pa_\J^\A}\longrightarrow B^2_{\prim}\overset{\ol\pa_\J^\A}\longrightarrow B^3\overset{\ol\pa_\J^\A}\longrightarrow \cdots \overset{\ol\pa_\J^\A}\longrightarrow B^{2n}\to 0,
$$
\bgn{proposition}
The complex $B^\bullet$ is the deformation complex of Einstein-Hermitian generalized connections and the infinitesimal space is given by the first cohomology of the complex $B^\bullet$.
\end{proposition}
\bgn{proof}The image of $d^\A(B^0)$ is the tangent space of the orbit of the action of the gauge group
$U(E,h)$ on the space $\GCon(E, h)^{hol}$. 
Let $\D_{\A_t}$ be a family of Einstein-Hermitian generalized connections in 
$\wtil {\M(E, h)}^{hol}_{EH}.$ Thus we have 
$$
\pi_{U^{-n}}^{\Herm}\F_{\A_t}(\psi)=\lam \id_E,
$$
where $\lam$ is a constant. 
Then we have 
$$
\pi_{U^{-n}}^{\Herm}d^\A ({\dot \A}\cdot\psi)=0.
$$  
Hence it follows that $K(\dot\A)=0.$
Since $\ol\pa_\J^{\A_t}\circ \ol\pa_\J^{\A_t}=0$, we also have 
$\ol\pa_\J^\A{\dot \A^{0,1}}=0$.
Thus the tangent space of $\wtil {\M(E, h)}^{hol}_{EH}$ at $\D^\A$ is given by 
the kernel of tha map $K\oplus \ol\pa_\J^\A: B^1\to B^2$
The moduli space ${\M(E, h)}^{hol}_{EH}$ is the quotient of  $\wtil {\M(E, h)}^{hol}_{EH}$ at $\D^\A$ by the action of the gauge group $U(E, h).$ Hence 
the tangent space of the moduli space is the quotient space 
$(\ker\,K\oplus \ol\pa_\J^\A: B^1\to B^2) /( \Im d^\A :B^0\to B^1).$
Hence we have the result.
\end{proof}
\bgn{proposition}
The complex $B^\bullet$ is an elliptic complex. 
Each cohomology group is finite dimensional over a compact manifold $M$.
\end{proposition}
\bgn{proof}To prove that the complex $B^\bullet$ is elliptic, we have to show that the sequence $B^\bullet$ is exact at the symbol level.
Let $\t\neq 0$ be a real cotangent vector at a point $x\in M.$ Then $\t$ is decomposed into $\t=\t^{1,0}+\t^{0,1},$ where 
$\t^{1,0}\in \L_\J$ and $\t^{0,1}=\ol{ \t^{1,0}}\in \ol \L_\J.$ 
By using a generalized K\"ahler structure $(\J, \J_\psi),$ we also have the decomposition
$$
\t=\t^{1,0}_+ + \t^{1,0}_-+\t^{0,1}_+ +\t^{0,1}_-,
$$
where $\t^{1,0}_\pm\in \L_\J^\pm, \,\, \t^{0,1}_\pm\in \ol\L_\J^\pm$ and 
$\t^{1,0}=\t^{1,0}_++\t^{1,0}_-$. 
Then the symbol $\sig$ is given by 
\bgn{align}
&\sig(d^\A, \t) \xi =\xi \otimes \t\\
&\sig(K\oplus \ol\pa_\J^\A, \t)\xi\otimes e =
\pi_{U^{-n}}^{\Herm}(\xi\otimes \t\cdot e\cdot \psi )\,\,\, +\xi\otimes(\t^{0,1}\w e^{0,1})\\
&\sig(\ol\pa_\J^\A, \t)(\zeta+\xi_2\otimes \a^{0,2} )=\xi_2\otimes(\t^{0,1}\w \a^{0,2})\\
&\sig(\ol\pa_\J^\A, \t)(\xi_i\otimes \a^{0,i})=\xi_i\otimes(\t^{0,1}\w\a^{0,i}),\qquad 
\text{ for }i\geq 3,
\end{align}
where $\xi\in u(E)_x$ and $\zeta \in \Herm(E, h)_x,$  and $ \xi_2,\xi_i (i\geq 3)\in\End(E)_x. $
The symbol complex is clearly exact at $B^0$ and $B^i(i\geq 3).$ 
To see that it is exact at $B^1$, we assume that $\sig(K\oplus \ol\pa_\J^\A, \t)a=0$ for $a\in B^1$, 
where $a=a^{1,0}+ a^{0,1}$ and $a^{0,1}\in \End(E)\otimes \ol\L_\J$. 
Then $\t^{0,1}\w a^{0,1}=0$. Thus there is a $f\in \End(E)$ such that $a^{0,1}=f\t^{0,1} .$
Then we also obtain $a^{1,0}=- f^*\t^{1,0}$ since $\tau (a)=-a$.
(Note that $a\in u(E)\otimes (\TT)$ satisfies $\tau (a) =\tau(a^{1,0})+\tau(a^{0,1})=-a.$ 
Thus $\tau (a^{1,0})=-a^{0,1}.$)  
Since the symbol of the map $\sig(K, \t)a=0$, we have $ \pi_{U^{-n}}^{\Herm}\t\w a\cdot\psi=0$.
Since $\t\in T^*_M,$ we have $\t\w\t=0.$ Thus $\t^{1,0}\w \t^{0,1}+\t^{0,1}\w \t^{1,0}=0.$ 
It follows that we have 
\bgn{align}
\sig(K, \t)a=&\pi_{U^{-n}}^{\Herm}\t\w a\cdot\psi\\
=&\pi_{U^{-n}}^{\Herm}\t\w(-f^*\t^{1,0} + f\t^{0,1})\cdot\psi=2 \pi_{U^{-n}}^{\Herm}((f+ f^*)\t^{1,0}\w\t^{0,1}\cdot\psi)=0
\end{align}
By using the decomposition $\t=\t^{1,0}_+ + \t^{1,0}_-+\t^{0,1}_+ +\t^{0,1}_-,$ we have 
 \bgn{align}
 \pi_{U^{-n}}\t^{1,0}\w\t^{0,1}\cdot\psi =&\t^{1,0}_+\cdot \t^{0,1}_+\cdot\psi\in U^{0,-n}\\
 =&-2\lan \t^{1,0}_+, \, \t^{0,1}_+\ran_{\tt}\cdot\psi\
 \end{align}
 where $\lan \t^{1,0}_+, \, \t^{0,1}_+\ran_{\tt}$ is real. 
Thus we obtain 
\bgn{align}
\sig(K, \t)a=&-4\lan \t^{1,0}_+, \, \t^{0,1}_+\ran_{\tt}\pi^{\Herm}  (f+ f^*)\cdot\psi\\
=&-4\lan \t^{1,0}_+, \, \t^{0,1}_+\ran_{\tt} (f+ f^*)\cdot\psi\\
=&0
\end{align}
Note that the projection $\pi^{\Herm}$ is the projection from $\End(E)$ to $\Herm(E,h)$. 
Since $\t\neq 0,$ we have $\lan \t^{1,0}_+, \, \t^{0,1}_+\ran_{\tt}\neq 0$. 
Thus $f+ f^*=0.$ 
Hence $f\in u(E)$. Then 
$a\in B^1$ is given by 
$$
a=f\t, \qquad f\in u(E)
$$
Thus the symbol complex is exact at $B^1.$
It is clear that the complex $(\w^\bullet \End(E)\otimes \ol \L_\J, \ol\pa_\J^\A)$ is elliptic. 
Thus the symbol complex of $(\End(E)\otimes\w^\bullet \ol \L_\J, \ol\pa_\J^\A)$ is exact and 
the alternating sum $\sum_i (-1)^i\dim \End(E)\otimes\w^i \ol \L_\J=0.$
Since $\End(E)=u(E)\oplus\Herm(E,h)$, it follows from the definition of $B^\bullet$ that we have 
\bgn{align}
&\dim B^0_x+\dim (B^2_{\prim})_x=\dim \End(E)_x+ \dim \End(E)\otimes\w^2\ol\L_{\J,x}\\
&\dim B^{1}_x=\dim \End(E)\otimes\ol\L_{\J,x}, \\
&\dim B^i_x=\End(E)\otimes\w^i\ol\L_{\J,x}\qquad (i\geq 3)
\end{align}
Then the exactness at $(B^2_{\prim})_x$ follows from the vanishing of the following alternating sum 
$$
\sum_{i=0}^{2n} (-1)^i \dim B^i_x =0.
$$
\end{proof}
\bgn{remark}
Let $C^{0,i}:=\End(E)\otimes \w^i\ol {\L_\J}$ be the Lie algebroid complex . 
Then there is a bundle map between $B^\bullet$ and $C^{0,\bullet}$
\bgn{equation*}
\xymatrix{
0\ar[r] &B^{0}\ar[d]\ar[r]^{d^{\D_\A}}&B^{1}\ar[d]\ar[r]^{K\oplus\ol\pa_\J^\A}&B^2\ar[d]\ar[r]^{\ol\pa_\J^\A}&B^3\ar[d]\ar[r]^{\ol\pa_\J^\A}&\cdots\\
0\ar[r] &C^{0,0}\ar[r]^{\ol\pa_\J^\A}&C^{0,1}\ar[r]^{\ol\pa_\J^\A}&C^{0,2}\ar[r]^{\ol\pa_\J^\A}
&C^{0,3}\ar[r]^{\ol\pa_\J^\A} &\cdots.}
\end{equation*}
\end{remark}
\section{Einstein-Hermitian co-Higgs bundles}
Let $(M,J, \ome)$ be a K\"ahler manifold and $(E, h)$ a Hermitian vector bundle over $M.$ 
We denote by $(\J_J, \J_\psi)$ the induced generalized K\"ahler structure, where $\psi=e^{\sqrt{-1}\ome}$. 
A generalized Hermitian connection $\D^{\A}=D^A+V$ consists of a Hermitian connection and a skew-symmetric endmorphism $V\in u(E).$ 
$\D^\A$ is decomposed as 
$$
\D^\A=(\D^\A)^{1,0}+(\D^\A)^{0,1}, 
$$
where $(\D^\A)^{1,0}: \Gam(E)\to \Gam(E\otimes \L_{\J_J})$ and $(\D^\A)^{0,1}: \Gam(E)\to \Gam(E\otimes \ol{\L_{\J_J}}).$
We assume that $(\D^\A)^{0,1}\circ (\D^\A)^{0,1}=0.$ Then $(\D^\A)^{0,1}$ gives a generalized holomorphic structure $\ol\pa_{\J_J}^\A$ on $E.$
A generalized holomorphic structure on $(M, \J_J, \J_\psi)$ is called {\it a co-Higgs bundle} \cite{Hi_2011-1}, \cite{Ra1_2013}, \cite{Ra2_2014}.
We shall describe the Einstein-Hermitian condition for a generalized Hermitian connection $\D^\A$ associated with 
a co-Higgs bundle on $E$. The Einstein -Hermitian condition
is given by 
$$
\pi^{\Herm}_{U^{-n}}\F_\A(\psi)=\id_E,
$$
where
$\F_\A(\psi)=F_A\cdot\psi+ d^A(V\cdot\psi)+\frac12[V\cdot V]\cdot\psi.$
Let $\Lam_\ome$ be the contraction by the K\"ahler form $\ome$. 
Since $(\Lam_\ome F_A)\ome^n=n F_A\w \ome^{n-1},$ 
then we have 
\bgn{align}
\pi^{\Herm}_{U^{-n}}F_A\cdot\psi=\frac{\lan F_A\cdot\psi,\, \ol\psi\ran_s }{\lan \psi, \, \ol\psi\ran_s}\psi
=-\frac{\sqrt{-1}}2\Lam_\ome F_A
\end{align}
$V$ is written as $V=\sum_i V_i v_i$, where $V_i \in u(E)$ and $v_i$ is a real vector field.
Then we have 
\bgn{align}
\pi_{U^{-n}}^{\Herm}d^A(V\cdot\psi)=&\pi^{\Herm}_{U^{-n}}(dV\cdot\psi +[A_\a\,\cdot \, V]\cdot\psi)\\
=&\pi^{\Herm}_{U^{-n}}\(\sum_i (d^A V_i) \sqrt{-1}i_{v_i}\ome \w\psi+\sqrt{-1}V_i\otimes (di_{v_i}\ome)\w\psi\)
\end{align}
Since $d^AV_i$ is a $u(E)$-valued $1$-form, we have $\pi_{U^{-n}}^{\Herm }(d^A V_i) \sqrt{-1}i_{v_i}\ome \w\psi=0.$
Hence we have $\pi_{U^{-n}}^{\Herm}d^A(V\cdot\psi)=0.$
$V$ is also written as $V=V^{1,0}+V^{0,1},$ where $V^{1,0}\in \End(E)\otimes\L_{\J_J}$ and 
$V^{0,1}\in \End(E)\otimes\ol{\L_{\J_J}}$ and $V^{1,0}=-(V^{0,1})^*$.
The K\"ahler form is locally written as 
$\sqrt{-1}\ome =\sum_i \t\w\ol{\t^i}$ for a basis $\{\t^i\}$ of $\w^{1,0}_M.$ 
We take the dual basis $\{z_i\}$ of $T^{1,0}_M$ of $\{\t^i\}$ and then $V^{1,0}$ is given by 
$V^{1,0}=\sum_i V^{1,0}_i z_i$ and $V^{0,1}=\sum_{\ol j}V^{0,1}_{\ol j}\ol z_j.$ 
Then we have 
\bgn{align}
\frac12[V\cdot V]\cdot\psi =\sum_i[ V^{1,0}_i, \, V^{0,1}_{\ol j}]z_i\w \ol{z_j}\psi
\end{align}
Since $\psi=e^{\sum_l \t^l\w\ol \t^l},$ we have 
$$
(z_i\w \ol{z_j})\cdot\psi=-z_i\cdot \t^j\cdot\psi=(-\del_{i,j}+\t^j\w\ol{\t^i})\cdot\psi.
$$
We also have 
\bgn{align}
\pi_{U^{-n}}(\t^i\w \ol{\t^j})\cdot\psi=&\frac{\lan (\t^i\w \ol{\t^j})\cdot\psi, \, \ol\psi\ran_s}{\lan\psi, \,\ol\psi\ran_s}
\psi\\
=&\frac{\lan (\t^i\w\ol{\t^i})\cdot\psi, \,\ol\psi\ran_s}{\lan\psi, \,\ol\psi\ran_s}\del_{i,j}\\
=&\frac{(\t^i\w\ol{\t^i})\w(\sum_l\t^l\w\ol{\t^l})^{n-1}2^{n-1}}{(n-1)! } 
\frac{n!}{2^n(\sum_l\t^l\w\ol{\t^l})^{n}}\del_{i,j}\\
=&\frac12\del_{i,j}
\end{align}
Thus we have 
$$
\pi_{U^{-n}}(z_i\w \ol{z_j})\cdot\psi=-\frac12 \psi
$$
Since $[ V^{1,0}_i, \, V^{0,1}_{\ol j}]=-[ V^{1,0}_i, \, (V^{1,0})^*]$ is Hermitian, we have 
$$
\pi_{U^{-n}}^{\Herm}\frac12[V\cdot V]\cdot\psi=-\frac12\sum [ V^{1,0}_i, \, V^{0,1}_{\ol i}]
$$
Hence we have 
\bgn{proposition}[Einstein-Hermitian conditions co-Higgs bundle]\label{Einstein-Hermitian conditions co-Higgs bundle}
The Einstein-Hermitian conditions for $(\J_J, \J_\psi)$ is given by 
$$
\sqrt{-1}\Lam_\ome F_A+\sum_i [V_i^{1,0}, \, (V_i^{1,0})^*]=\lam \id_E,
$$
where $\lam$ is a constant.
\end{proposition}
It is remarkable that the Einstein-Hermitian condition coincide with the equation in \cite{Hi_2011-1} which is an analogy of the one of Higgs bundles
. 
However if we take $\psi=e^{b+\sqrt{-1}{\ome}}$ for a $d$-closed real $2$-form $b$, 
then the Einstein-Hermitian condition is different from the equation in \cite{Hi_2011-1}. 
In fact, the term $d^A( V\cdot\psi)$ is given by 
$$
d^A (V\cdot\psi)=\sum_j d^A(V_j\otimes (i_{v_i}b+\sqrt{-1}i_{v_j}\ome))\cdot\psi
$$
Thus we have 
$$
\pi_{U^{-n}}^{\Herm}d^A(V\cdot\psi) =\pi_{U^{-n}}\sum_j d^A(V_j \otimes i_{v_j}b)\cdot\psi
$$
This term does not vanish in general. 
The third term $\pi_{U^{-n}}^{\Herm}[V\cdot V]\cdot\psi$ changes by the action of $b$-field. 
As in Proposition \ref{prop: b-field action} and Proposition \ref{prop: invariance of b-field action}, our Einstein-Hermitian condition is equivalent under the action of $b$-field. 
By using an action of $b$ field, we shall describe our Einstein-Hermitian condition in the cases of $\psi=e^{b+\sqrt{-1}\ome}.$ 
An Einstein-Hermitian generalized connection $\D^\A$ is changed to 
$\Ad_{e^{-b}}\D^\A=(d+A - \ad_b V)+ V$ which is an Einstein-Hermitian co-Higgs bundle over $(M, e^{\sqrt{-1}\ome}).$
From Proposition \ref{Einstein-Hermitian conditions co-Higgs bundle}, 
our Einstein-Hermitian condition is the following 
\bgn{equation}\label{eq:nonabelian Kahler Ricci soliton}
\sqrt{-1}\Lam_\ome \(F_A-d^A(\ad_b V)+\frac12 [\ad_b V, \, \ad_b V]\)+\sum_i [V_i^{1,0}, \, (V_i^{1,0})^*]=\lam \id_E,
\end{equation}
In the case of $b=\ome$, the equation 
(\ref{eq:nonabelian Kahler Ricci soliton}) is regarded as a $\End(E)$-valued K\"ahler Ricci Solitons.

\medskip
\noindent
E-mail address: goto@math.sci.osaka-u.ac.jp\\
\noindent
Department of Mathematics, Graduate School of Science,\\
\noindent Osaka University Toyonaka, Osaka 560-0043, JAPAN 
\end{document}